\documentclass{amsart}

\usepackage{amssymb}
\usepackage{eucal}
\usepackage{amsfonts}
\usepackage[frame,ps,dvips,matrix,arrow,curve,rotate]{xy}

\newtheorem{thm}{Theorem}[section]
\newtheorem{cor}[thm]{Corollary}
\newtheorem{lem}[thm]{Lemma}
\newtheorem{prop}[thm]{Proposition}
\theoremstyle{definition}
\newtheorem{defn}[thm]{Definition}
\theoremstyle{remark}
\newtheorem{rem}[thm]{Remark}
\theoremstyle{remark}
\newtheorem{ex}[thm]{Example}

\newcommand{\X}{\mathcal{X}}
\newcommand{\Y}{\mathcal{Y}}
\newcommand{\Z}{\mathcal{Z}}

\newcommand{\Ss}{\mathcal{S}}

\newcommand{\Xm}{X_{mod}}
\newcommand{\Ym}{Y_{mod}}

\newcommand{\ox}{\omega_x}
\newcommand{\oy}{\omega_y}
\newcommand{\OX}{\Omega_{\X}}

\newcommand{\pim}{\pi_{mod}}

\newcommand{\Ra}{\Rightarrow}

\newcommand{\rsa}{\rightsquigarrow}

\newcommand{\x}{\times}
\newcommand{\st}[1]{\stackrel{#1}{\rightrightarrows}}

\newcommand{\PX}{\Pi_1(\X)}
\newcommand{\pX}{\pi_1(\X,x)}
\newcommand{\PhX}{\Pi_1^{h}(\X)}
\newcommand{\phX}{\pi_1^{h}(\X,x)}

\newcommand{\pXm}{\pi_1(\Xm,x)}

\newcommand{\PY}{\Pi_1(\Y)}
\newcommand{\pY}{\pi_1(\Y,y)}
\newcommand{\PhY}{\Pi_1^{h}(\Y)}
\newcommand{\phY}{\pi_1^{h}(\Y,y)}

\newcommand{\ppXX}{\pi_1(X,x)}
\newcommand{\ppYY}{\pi_1(Y,y)}

\newcommand{\pXX}{\pi_1(X',x')}

\def\smashedlongrightarrow{\setbox0=\hbox{$\longrightarrow$}\ht0=1pt\box0}
\def\risom{\buildrel\sim\over{\smashedlongrightarrow}}
\def\smashedst{\setbox0=\hbox{$\rightrightarrows$}\ht0=4pt\box0}
\newcommand{\sst}[1]{\stackrel{#1}{\smashedst}}

\let\:=\colon

\newcommand{\Spec}{\operatorname{Spec}}
\newcommand{\Gal}{\operatorname{Gal}}

\newcommand{\Ob}{\operatorname{Ob}}
\newcommand{\Mor}{\operatorname{Mor}}

\newcommand{\Ker}{\operatorname{ker}}

\newcommand{\im}{\operatorname{im}}

\begin{document}

\title[Fundamental Groups of Algebraic Stacks]{Fundamental Groups of Algebraic Stacks}%
\author{B. Noohi}%


\begin{abstract}
We study fundamental groups of algebraic stacks. We show that
these fundamental groups carry an additional structure coming from
the inertia groups. Then use this additional structure to analyze
geometric/ topological properties of stacks. We give an explicit
formula for the fundamental group of the coarse moduli space.
Also, we use these additional structures to give a necessary and
sufficient for an algebraic stack to be uniformizable.

\end{abstract}
\maketitle
\section{Introduction}

In \cite{SGA1} Grothendieck introduces a general formalism to
associate a {\em fundamental group} with a pointed scheme $X$.
This is a profinite group which should be thought of as the
``profinite completion of the (virtual) Poincar\'{e} fundamental
group'' of $X$. Grothendieck's ideas are indeed general enough to
apply to  any reasonable notion of ``space'' (e.g., a topos). For
instance, back to topology, one can apply the theory to (connected
locally 1-connected) topological spaces and what comes out is the
profinite completion of the usual fundamental group.
Grothendieck's theory also applies to Algebraic stacks, for
instance via the topos associated to the \'{e}tale or smooth site
of the stack \footnote{Finally there is a written account of the
topos theoretic approach. See the articles by V. Zoonekynd in
arXiv:math.AG.}.  The aim of this paper is to write out the theory
in a rather detailed way, but without using the language of topoi,
and then to explore the features that are special to this
particular case.

Let $\X$ be a connected Algebraic stack. Let $x$ be a geometric
point of $\X$, and let $\rho_x$ be its residue gerbe.
 The following simple observation is essential:
 {\em The automorphism group of the gerbe $\rho_x$ is isomorphic
 to the geometric fundamental group
of $\rho_x$; in particular, we have a natural group homomorphism
$\ox$ from the automorphism group of $\rho_x$ to the fundamental
group of $\X$ at $x$.} This means that the fundamental group of an
algebraic stack at a point $x$ comes equipped with an extra
structure -- namely, the map $\ox$. The significance of this map
is that it relates local data (the automorphism group of the
residue gerbe) to global data (the fundamental group). Since the
fundamental groups of $\X$ at different points are isomorphic
(with an isomorphism that is unique up to conjugation), we can
indeed map the automorphism group of the residue gerbe of any
point $x'$ into $\pX$. The normal subgroup generated by all the
images will be a well-defined subgroup of $\pX$. Let $N$ be the
closure of this subgroup. We prove the following

\begin{thm}
The group $\pX/N$ is naturally isomorphic to the fundamental group
of the moduli space of $\X$ (see Theorem \ref{G5} for the precise
statement).
\end{thm}

The maps $\ox$ can be used to detect whether a given algebraic
stack is uniformizable:

\begin{thm}
A Deligne-Mumford stack $\X$ is uniformizable if and only if all
the maps $\ox$ are injective (see Theorem \ref{T:uniformization}
for the precise statement).
\end{thm}

An algebraic stack being {\em uniformizable} means that it has a
finite \'{e}tale representable cover by an algebraic space
(roughly speaking, its ``universal cover'' is an algebraic space).

This paper is organized into two parts. Part one consists of the
main construction and the main results. Part two should be thought
of as a technical companion to part one where we   strengthen the
results of part one by introducing some more elaborate techniques.

\tableofcontents

\section{Review and Conventions}{\label{S:review}}

Our reference for the theory of stacks is \cite{LM}. We quickly
review a few basic facts that we will need, some of which are not
explicitly mentioned in \cite{LM}.

We begin with some conventions:

In this paper, whenever we use the word {\it category} for
something that is really a 2-category, we simply mean the category
that is obtained by identifying 2-isomorphic 1-morphisms. Typical
example: the 2-category of algebraic stacks (say, over a fixed
base algebraic space $S$).

We use calligraphic symbols ($\X$, $\Y$,...) for algebraic stacks
and ordinary symbols  ($X$, $Y$,...) for algebraic spaces or
schemes. We denote the Zariski topological space of an algebraic
stack $\X$ by $|\X|$.

Throughout this paper the word {\em representable} means
representable by algebraic spaces. Although it is an extraneous
condition in most part of the paper, but we will assume that our
algebraic stacks are locally Noetherian.

By a smooth (respectively, flat) {\em  chart} for an algebraic
stack $\X$ we mean a smooth (respectively, flat and of finite
presentation) surjective map $p \: X \to \X$, where $X$ is an
algebraic space. When we do not specify an adjective (flat or
smooth) for a chart, we mean a smooth chart. To a smooth
(respectively, flat) chart, we can associate a smooth
(respectively, flat) groupoid $R \sst{s,t} X$, where $R:=X
\x_{\X}X$, and the maps $s,t$ are the two projections. Conversely,
given a  flat groupoid $R:=X \x_{\X}X$ on an algebraic space $X$,
we can construct a quotient stack $\X=[X/R]$ (see Artin's theorem
in \cite{LM}, Section 10). These two constructions are inverse to
each other.

A useful fact to keep in mind is the following: Let $X \to \X$ be
a chart for $\X$, and let $R_X \st{} X$ be the corresponding
groupoid. Let $\Y \to \X$ be a representable morphism, and let $Y
\to \Y$ be the pull-back chart for $\Y$. Let $R_Y \st{} Y$ be the
corresponding groupoid. Then, the  diagram
\begin{equation}{\label{0}}
\xymatrix{ R_Y \ar[r] \ar[d] & R_X \ar[d] \\
             Y \ar[r]        & X }
\end{equation}
\noindent is cartesian, where the vertical arrows are either
source or target maps.

A {\em covering space} for an algebraic stack $\X$ is a finite
\'{e}tale representable morphism $f \: \Y \to \X$, where $\Y$ is a
connected algebraic stack. We sometimes call such a map  a {\em
covering map}.

We say a morphism $f \: \X \to \Y$ of algebraic stacks is a {\em
monomorphism}, if for any scheme $T$, the induced map $\X(T) \to
\Y(T)$ on the groupoids of $T$-points is fully faithful. When $f$
is representable, any base change of $f$ to a map of schemes will
be a monomorphism in the usual sense, and vice versa. In
particular, a locally closed immersion is a monomorphism.

\subsection{Stabilizer group of an algebraic stack}{\label{S:Stab}}

Let $R \st{} X$ be a flat groupoid over an algebraic space $X$,
and let $\X=[X/R]$ be the the corresponding quotient stack. The
{\em stabilizer group} of $R \st{} X$, or simply of $\X$ (a
harmless abuse of terminology!), is the group space  $S \to X$
that is defined by the following cartesian diagram:

$$\xymatrix{ S_X \ar[r] \ar[d] & R \ar[d]^{(s,t)} \\
             X \ar[r]_(0.35)\Delta & X \x X}$$

\noindent or, equivalently, by the following cartesian diagram

$$\xymatrix{S_X \ar[r] \ar[d] & \mathcal{S}_{\X} \ar[r] \ar[d] &
\X \ar[d]^{\Delta} \\
          X \ar[r]_p & \X \ar[r]_(.34){\Delta} & \X \x \X}$$

\noindent The group stack $\mathcal{S}_{\X} \to \X$ is called the
{\em inertia stack} or the {\em stabilizer group stack} or the
{\em automorphism group stack} of $\X$. The group stack
$\mathcal{S}_{\X} \to \X$ is representable. In particular, any
property of (representable) morphisms of stacks can be attributed
to the stabilizer of an algebraic stack (e.g., properties such as
finite, quasi-finite, \'{e}tale, unramified, reduced/connected
geometric fibers, and so on). The stabilizer group stack of an
algebraic stack $\X$ is always of finite type and separated;  when
$\X$ is a Deligne-Mumford stack it is unramified (hence
quasi-finite) as well.  When $x \: \Spec k \to \X$ is a geometric
point, then the {\em stabilizer group scheme of $x$}, denoted
$S_x$, is defined to be the fiber of $\X \to \X \x \X$ over the
diagonal point $x \to \X \x \X$ (equivalently, $S_x$ is the fiber
of $\mathcal{S}_{\X} \to \X$ over $x$). If $x_0$ is a lift of $x$
to a chart $X \to \X$ for $\X$, then $S_x$ is naturally isomorphic
to the fiber of $S_X \to X$ over $x_0$. We can also talk about the
stabilizer group of a point in the underlying space $|\X|$, but it
will be defined only up to isomorphism (unless we fix a geometric
point for it).

Throughout the text, whenever we use the phrase `stabilizer group
of $\X$', the reader can either think of the  stabilizer group
stack $\Ss_{\X} \to \X$, defined as above, or,  those who do not
like the word ``group stack'', can implicitly  fix a chart and
work with the stabilizer group of the corresponding groupoid.

The following two lemmas make it easy to play around with
stabilizer groups. Proofs are  easy consequences of definitions
and are left to the reader.

\begin{lem}{\label{P4.7}}
Let $\X$ be an algebraic stack, and let $p \: X \to \X$ be a chart
for it. Let $X' \to \X$ be another  chart for $\X$ that factors
through $X$. Then, we have the following cartesian diagram of
stabilizer group spaces:

$$\xymatrix{ S_{X'} \ar[r] \ar[d] & S_X \ar[d] \\
             X' \ar[r] & X}$$
\end{lem}

Let $\X$ be an algebraic stack and let $p \: X \to \X$ be a chart
for it. Let $f \: \Y \to \X$ be a representable morphism of
stacks, and let $q \: Y \to \Y$ be the chart for $\Y$ obtained by
pulling back $p$ via $f$. There is a natural group homomorphism
$S_Y \to  S_X \x_X Y$ (as group spaces over $Y$). The following
lemma says that, as far as algebro-geometric properties of
morphisms are concerned, this map behaves like the diagonal map
$\Y \to \Y \x_{\X} \Y$.

\begin{lem}{\label{P4.8}}
There is a natural cartesian diagram

$$\xymatrix{ S_Y \ar[r] \ar[d] &  S_X \x_X Y \ar[d] \\
              \Y \ar[r]_(0.4){\Delta} & \Y \x_{\X} \Y } $$
\end{lem}

\begin{cor}{\label{P4.9}}

Let $f \: \Y \to \X$ be a morphism of stacks. Let $f \: \Y \to \X$
be a representable morphism of stacks, and let $q \: Y \to \Y$ be
the chart for $\Y$ obtained by pulling back $p$ via $f$. Then:

\begin{itemize}
\item[$\mathbf{i})$]  If $f$ is unramified, then $S_Y$ is naturally
isomorphic to an open subgroup space of the pull back group space
$S_X \x_X Y$ (as group spaces over $Y$). In particular, the
stabilizer group of a geometric point $y$ of $\Y$ is isomorphic to
an open subgroup of that of $f(y)$. In particular, the stabilizer
group of $y$ is reduced if and only if the stabilizer group of
$f(y)$ is so. Finally, $\Y$ has reduced stabilizers if and only if
$\X$ does.

\item[$\mathbf{ii})$] If $f$ is separated, then $S_Y$ is naturally
isomorphic to a closed subgroup space of the pull back group space
$S_X \x_X Y$ (as group spaces over $Y$). In particular, if $\X$
has finite (respectively, proper) stabilizer, then so does $\Y$.
\end{itemize}
\end{cor}

\subsection{Groupoids}{\label{S:groupoid}}

The language of groupoids turns out to the most natural one to
formulate our results. In this section we fix some notations and
terminology related to groupoids.

Let $\Pi$ be a groupoid. We use the words {\em object} and {\em
point} interchangeably to refer to objects of $\Pi$ (viewed as a
category). Similarly, we use the words {\em morphism} and {\em
arrow} interchangeably to refer to morphisms of $\Pi$. For a point
$x \in \Pi$, we denote the automorphism group of $\Pi$ at $x$ by
$\Pi(x)$.

Let $\Pi$ be a groupoid, and let $\Pi'$ be a subgroupoid of it. We
say that $\Pi'$ is a {\sl normal} subgroupoid, if for every  arrow
$b$ in $\Pi'$ and every arrow $a$ in $\Pi$, we have $aba^{-1} \in
\Pi'$, whenever the composition is defined. For an arbitrary
subgroupoid $\Pi' \subseteq \Pi$, we define its {\sl normal
closure}, denoted $N\Pi'$, to be the smallest normal subgroup of
$\Pi$ containing $\Pi'$.

The above definition of a normal subgroup is equivalent to the
following: For every pair of points $x_1$ and $x_2$ in $\Pi'$ and
every arrow between them, the induced isomorphism $\Pi(x_1) \to
\Pi(x_2)$ maps $\Pi'(x_1)$ isomorphically to $\Pi'(x_2)$.

\pagebreak

{\LARGE\part{Basic results}}{\label{C:basic}}

\vspace{0.5in}

\section{Hidden Paths in Algebraic Stacks}

Algebraic stacks form a 2-category, meaning that, given a pair of
morphism with the same source and target, we could talk about
transformations (or 2-morphisms) between them. These 2-morphisms
behave somewhat like {\em homotopies}. In this section we will
exploit this idea systematically. In the next section, we see how
to obtain actual homotopies out of these transformations. In this
paper we are only interested in the case where the source is the
spectrum of an algebraically closed field. We make the following

\begin{defn}{\label{D:hidden}}
Let $x, x' \: \Spec k \to \X$ be two geometric points. By a {\sl
hidden path} from $x$ to $x'$, denoted $x \rsa x'$, we mean a
transformation from $x$ to $x'$. The {\sl hidden fundamental
group} of $\X$ at $x$, denoted $\phX$, is the group of
self-transformations of $x$. We define the {\sl hidden fundamental
groupoid} of $\X$, denoted $\PhX$, as follows:

\begin{itemize}
\item $\Ob \PhX$=\{geometric points of $\X$\}
\item $\Mor(x,x')$=\{hidden paths $x \rsa x'$\}.
\end{itemize}
\end{defn}

The multiplication $\gamma_1 \gamma_2$ of hidden paths $\gamma_1
\in \Mor(x,x')$ and $\gamma_2 \in \Mor(x',x'')$ is defined by
composition of transformations. For any geometric point $x \:
\Spec k \to \X$, the group of automorphism of $x$, viewed as an
object of the category $\PhX$, is equal to the hidden fundamental
group $\phX$.

The following proposition follows immediately from the definition.

\begin{prop}{\label{P:mono}}
A monomorphism (see end of Section \ref{S:review}) of algebraic
stacks induces an isomorphism of hidden fundamental groups.
\end{prop}

\begin{defn}{\label{D:pointed}}
By a {\sl pointed map} (or {\em a pointed morphism}) $(f, \phi) \:
(\Y,y) \to (\X,x)$ we mean a pair $(f, \phi)$, where $f \: \Y \to
\X$ is a morphism of algebraic stacks and $\phi \: x \rsa f(y)$ is
a hidden path (Definition \ref{D:hidden}).
\end{defn}

For a morphism $f \: \Y \to \X$ of algebraic stacks, we obtain an
induced  map of groupoids $\Pi_1^h(f) \: \PhY \to \PhX$. So, for a
pointed map $(f, \phi) \: (\Y,y) \to (\X,x)$, we obtain a natural
map $\pi_1^{h}(f, \phi)$ between the hidden fundamental groups as
follows:
\begin{align*}
\pi_1^{h}(f, \phi) \: \pi_1^{h}(\Y, y) & \to \pi_1^{h}(\X, x) \\
   \beta \ \ \  & \mapsto  \phi f(\beta) \phi^{-1}.
\end{align*}

The following proposition gives a simple description of the hidden
fundamental groups.

\begin{prop}{\label{P4}}
Let $p \: X \to \X$ be a chart for $\X$, and let $x \: \Spec k \to
\X$ be a geometric point  of $\X$. Let $x_0$ be a geometric point
of $X$, and let $\phi \: x \rsa p(x_0)$ be a hidden path (i.e.,
$(p,\phi) \: (X,x_0) \to (\X,x)$ is a pointed map) . Then we have
a natural isomorphism $S_{x_0} \risom \phX$, where $S_{x_0}$
denotes the set of $k$-points of the fiber of the stabilizer group
space $S \to X$ at $x_0$. In particular, $\phX$ has a natural
structure of an algebraic group.
\end{prop}

\begin{proof}
By definition of the fiber product in the 2-category of algebraic
stacks, a $(\Spec k)$-point of $S_{x_0}$ corresponds to a
transformation of functors from $p \circ x_0$ to itself. So, we
have a natural isomorphism $\lambda \: S_{x_0} \to \pi_1^{h}(\X,
p(x_0))$. We define $S_{x_0} \risom \phX$ by sending $s \in
S_{x_0}$ to $\phi \lambda(s) \phi^{-1}$.
\end{proof}

\begin{rem}
As we saw in the proof, this isomorphism does depend on the choice
of $\phi$. Indeed, if $\phi' \: x \rsa p(x_0)$ is another hidden
path, then the two isomorphisms will be conjugate by the element
$\gamma=\phi' \phi^{-1} \in \phX$, that is, $\lambda_{\phi'}
=\gamma \lambda_{\phi} \gamma^{-1}$.
\end{rem}

By the above proposition, the hidden fundamental groupoid at a
geometric point $x$ is isomorphic to the automorphism group of the
residue gerbe at $x$. If $\X$ is defined as a quotient of a group
action on an algebraic space $X$ , then the hidden fundamental
group at  the point $x$ is isomorphic to the isotropy subgroup of
any point in $X$ lying above $x$. If $\X$ is defined by as the
moduli space of certain objects, then the hidden fundamental group
at a point $x$ is naturally isomorphic to the automorphism group
of the object represented by $x$.

\begin{defn}
We say that a geometric point $x$ of $\X$  (or its image in the
underlying space of $\X$) is {\sl unramified}, if $\phX$ is
trivial. We say it is {\sl schematic}, if the residue gerbe at $x$
is a scheme (necessarily of the form $\Spec k$).
\end{defn}

A point is schematic if and only if its stabilizer group is
trivial. It is unramified if and only if its stabilizer group is
connected and zero-dimensional. A point is schematic if and only
if it is unramified and its stabilizer is reduced. In particular,
these notions are equivalent for Deligne-Mumford stacks, or for
algebraic stacks over a field of characteristic zero. However, in
general they are not equivalent (Example \ref{P4.6} below).

It is worthwhile to keep in mind that a point $x$ being schematic
or unramified only depends on the residue gerbe at $x$.

\begin{cor}{\label{P4.5}}
Let $\X$ be an algebraic stack.
\begin{itemize}
\item[$\mathbf {i})$] If $\X$ has quasi-finite stabilizer,
then $\phX$ is finite for every geometric point $x$.

\item[$\mathbf {ii})$] $\X$ is an algebraic space if and only if
all its points are schematic.
\end{itemize}
\end{cor}

\begin{proof}
Immediate.
\end{proof}

\begin{ex}{\label{P4.6}}{\it An algebraic stack all whose points
are unramified need not be an algebraic space.} For instance, take
the classifying stack of a (non-trivial) connected
zero-dimensional group scheme over a field (necessarily of
positive characteristic). More explicitly, take an elliptic curve
defined over a field of positive characteristic and let it act on
itself via Frobenius. Then the quotient stack is an algebraic
stack whose underlying set is just a single point and the
stabilizer group at this point is isomorphic to the kernel of
Frobenius which is a connected finite flat group scheme supported
at a single point. In particular, $\phX$ is trivial by Proposition
\ref{P4}.
\end{ex}

The following lemma is straightforward.

\begin{lem}{\label{D1}}
Let $f \: \X \to \Y$ be a morphism of algebraic stacks. If $f$ is
representable, then for every geometric point $x$ of $\X$ the
induced map $\phX \to \pi_1^h(\Y,f(x))$ is injective.
\end{lem}

The converse of the this Proposition is true in zero
characteristic, or when $\X$ is a Deligne-Mumford stack. For a
counterexample, take the algebraic stack of Example \ref{P4.6} and
map it to the ground field.

\section{The Galois Category of an Algebraic Stack}{\label{S:Galois}}

In this section, we use Grothendieck's formalism of Galois
categories to associate with a pointed connected algebraic stack
$(\X,x)$ a fundamental group $\pX$. We refer the reader to
\cite{SGA1} for an account of the theory of Galois categories and
fundamental functors.

The reader who is familiar with the language of topos theory will
immediately realize that the Galois category that we associate to
a pointed algebraic stack is nothing but the Galois category of
locally constant sheaves in the pointed topos associated to the
stack (the choice of topology turns out to be immaterial).
Therefore, we obtain the same fundamental groups as we would have
obtained via the topos theoretic approach. The advantage of our
more explicit approach, is that it makes it clear how the
ramification of $\X$ at $x$ induces an extra structure on $\pX$.
More precisely, we will show that there is a natural group
homomorphism $\ox \: \phX \to \pX$. The significance of this map
is two-fold: on the one hand, it relates local data to global
data, and, on the other hand, it gives a homotopy theoretic
meaning to the ramification structure of stacks.

\begin{defn}{\label{D:covering}}
We say that a morphism $f \: \Y \to \X$ of connected algebraic
stacks is a {\em covering map}, or $\Y$ is a covering space of
$\X$ (via $f$), if $f$ is finite \'{e}tale and representable.
\end{defn}

Let $\X$ be a connected algebraic stack and let $x \: \Spec k \to
\X$ be a geometric point. We define the Galois category
$\mathbf{C}_{\X}$ of $\X$ and the fundamental functor $F_x$
associated to the point $x$ as follows:

\begin{itemize}
\item  $\Ob (\mathbf{C}_{\X}) = \left\{ \begin{array}{rcl}
         (\Y, f)   & \vert & \text{ $\Y$ \ an  algebraic stack;} \\
            & &  f \: \Y \to \X \ \  \text{a covering space.}
                                      \end{array} \right\}$

\item  $\Mor_{\mathbf{C}_{\X}}((\Y,f), (\Z,g))=\left\{\begin{array}{rcl}
           (a,\Phi) & \vert & a \: \Y \to \Z \  \ \text{a morphism;} \\
        && \Phi \: f \Ra g \circ a \  \ \text{a 2-morphism.} \end{array}\right\}_{/_\sim}$
\end{itemize}

\noindent where $\sim$ is defined by

$$(a, \Phi) \sim (b, \Psi) \ \ \text{if}\ \  \exists \Gamma  \: a
\Ra b \ \ \text{such that} \  g(\Gamma) \circ \Phi = \Psi.$$

\vspace{.2in}

\noindent The fundamental functor $F_x \: \mathbf{C}_{\X} \to
((Sets))$ is defined as follows

$$F_x(\Y)= \left\{\begin{array}{rcl} (y, \phi) & \vert & y \:
\Spec k \to \Y
                                                    \ \ \text{geometric point} \\
  & & \phi \: x \rsa f(y) \ \ \text{hidden path} \end{array} \right\}_{/_\sim}$$

\noindent where $\sim$ is defined by

$$(y, \phi) \sim (y',\phi') \ \ \text{if}\ \  \exists \beta \: y
\rsa y' \ \ \text{such that} \ f(\beta) \circ \phi =\phi'.$$

\noindent Note that the set $F_x(\Y)$ is in a natural bijection
with the underlying set of the geometric fiber $\Spec k
\times_{\X} \Y$, which is isomorphic to a disjoint union of copies
of $\Spec k$ .

Proof that  $\mathbf{C}_{\X}$ is a Galois category, with  $F_x$  a
fundamental functor for it, is straightforward (for a detailed
proof see \cite{Noohi}). Following \cite{SGA1}, we define $\pX$ to
be the group of self-transformations of the functor $F$. The group
$\pX$ is a profinite group. Recall that, by definition, a
profinite group is a compact Hausdorff totally disconnected
topological group. Equivalently, a profinite group is a group that
is isomorphic to a directed inverse limit of finite groups. There
is an equivalence of categories between $\mathbf{C}$ and the
category of $\pX$-sets, under which $F$ correspond to the
forgetful functor. If $x'$ is another geometric point, then
$\pi_1(\X,x')$ is isomorphic to $\pX$ via an isomorphism that is
unique up to conjugation by an element of $\pi_1(\X,x')$. The
fundamental group $\pX$ defined above classifies pointed covering
spaces of $(\X,x)$, in the sense that, there is a one-to-one
correspondence between isomorphism classes of pointed connected
covering spaces of $(\X,x)$ and open subgroups of $\pX$.

We define the fundamental groupoid $\PX$ of $\X$ as follows:

\begin{itemize}
\item $\Ob \PX$=\{geometric points of $\X$\}

\item $\Mor(x,x')$=\{transformations of functors $F_x \to F_{x'}$\}
\end{itemize}

\noindent The multiplication $\gamma_1 \gamma_2$ of ``paths''
$\gamma_1 \in \Mor(x,x')$ and $\gamma_2 \in \Mor(x',x'')$ is
defined to be the composition of transformations. The fundamental
groupoid is a connected groupoid whose group of automorphisms at
any point $x \in \PX$ is $\pX$.

Let $f \: \Y \to \X$ be an arbitrary covering space, and let $x \:
\Spec k \to \X$ and $x' \: \Spec k \to \X$ be two geometric
points. Let $\gamma \: x \rsa x'$ be a hidden path. Then we obtain
a map of sets $F_x(\Y) \to F_{x'}(\Y)$ by sending $(y, \phi)$ to
$(y, \gamma^{-1}\phi)$. This map is functorial in $\Y$; so it
gives rise to a transformation of functors $F_x \to F_{x'}$, that
is, an element in $\Mor_{\PX}(x,x')$. This construction produces a
natural  map of groupoids $\Omega \: \PhX \to \PX$. In particular,
for any geometric point $x$, we have a natural group homomorphism
$\omega_x \: \phX \to \pX$. These constructions are all functorial
in the following sense. Let $f \: \Y \to \X$ be a morphism of
algebraic stacks. Then we have natural maps of groupoids $\Pi_1(f)
\: \PY \to \PX$ and $\Pi_1^{h}(f) \: \PhY \to \PhX$ that are
compatible with $\Omega$, i.e., $\Pi_1(f) \circ \Omega_{\Y} =
\Omega_{\X} \circ \Pi_1^{h}(f)$. Similarly, if $(f, \phi) \:
(\Y,y) \to (\X ,x)$ is a pointed map, we obtain natural maps
$\pi_1(f, \phi) \: \pY \to \pX$  and $\pi_1^{h}(f, \phi) \: \phY
\to \phX$. More explicitly, for $\beta \in \pY$, $\pi_1(f,
\phi)(\beta)$ is defined to be $\phi \pi_1(f)(\beta) \phi^{-1}$.
The map $\pi_1^{h}(f, \phi) \: \phY \to \phX$, whose definition is
similar to that of  $\pi_1(f, \phi)$, has already been introduced
in the previous section. Once again we have the compatibility
relation $\pi_1(f, \phi) \circ \omega_y = \omega_x \circ
\pi_1^{h}(f, \phi)$.

\begin{ex}{\label{P5.9}}
Let $G$ be an algebraic group of finite type over an algebraically
closed field $k$. Let $X=\Spec k$, and let $\X=[X/G]$ be the
quotient of $X$ under the trivial action of $G$ (i.e., the
classifying stack of $G$) made into a pointed stack via the
quotient map $\Spec k \to \X$. Then, $\pX$ is naturally isomorphic
to the group $G/{G^0}$ of connected components of $G$, and $\phX$
is naturally isomorphic to $G$ itself (more precisely, the group
of $k$-points of $G$). The map $\ox$ is simply the quotient map $G
\to G/{G^0}$.  Proof that $\pX \cong G/{G^0}$ is easy. For
instance, we could use a fiber homotopy exact sequence argument
(see the appendix)
 to prove
the result. We give a more direct proof: We have a  Galois
covering $[\Spec k/G^{0}] \to \X$  whose Galois group is
$G/{G^0}$. The claim follows if we show that $[\Spec k/G^{0}]$ is
simply connected. So we may assume $G$ is connected. Let $\Y \to
\X$ be a covering map,  $Y \to \Y$ the pull-back of $X \to \X$,
and $R_Y \st{} Y$ the corresponding groupoid. Then, $Y$ is a
disjoint union of copies of $X$. Since $G$ is connected, it
follows from Corollary \ref{P4.9}, that the diagram $$\xymatrix{
S_Y \ar[r] \ar[d] & S_X \ar[d] \\
             Y \ar[r]        & X }$$
\noindent is cartesian. It now follows from the cartesian diagram
(\ref{0}) of Section \ref{S:review}, that $S_Y$ is equal to $R_Y$.
Therefore, $\Y$ is just a disjoint union of copies of $\X$, hence
the claim. (Note that this argument is valid whenever $G$ is a
connected, and $X$ is simply connected.)
\end{ex}

\section{Basic properties of hidden fundamental groups}{\label{S:basic}}

In this section we prove a few basic facts about hidden
fundamental groups. The following lemma is essential.

\begin{lem}{\label{P6}}
Let $f \: \Y \to \X$ be a covering  space. Let $y$ be a geometric
point of $\Y$, and let $x=f(y)$ be its image in $\X$. Then the
following diagram is cartesian:

$$\xymatrix{ \phY \ar[r]^{\omega_y} \ar[d]_{\pi_1^{h}(f)} & \pY
\ar[d]^{\pi_1(f)}\\
             \phX \ar[r]_{\omega_x} & \pX } $$
\end{lem}

\begin{proof}
Let $\gamma$ be in $\phX$. We want to show that, if
$\omega_x(\gamma)$ is in the image of $\pi_1(f)$, then there
exists a unique $\alpha \in \phY$ that maps to $\gamma$ via
$\pi_1^{h}(f)$. The uniqueness is obvious, since $\pi_1^{h}(f)$ is
injective by Lemma \ref{D1}. On the other hand, $\omega_x(\gamma)$
being in the image of $\pi_1(f)$ exactly means that, under the
action of $\gamma$ on $F_x(\Y)$, the point $(y, id)$ remains
invariant. That means, $(y,id) \sim (y, \gamma^{-1})$. Therefore,
there exists $\beta \in \phY$ such that $f(\beta)=\gamma^{-1}$.
The element $\alpha = \beta^{-1}$ has the desired property.
\end{proof}

\begin{rem}
The above lemma fails when $f$ is not representable. For instance,
let $\Y$ be the classifying space of an arbitrary non-trivial
finite group, and let $f$ be the moduli map. Then $\omega_x$ is
zero, but $\omega_y$ is not.
\end{rem}

\begin{cor}{\label{P5}}
The kernel of $\omega_x$ is equal to $\bigcap \im(\pi_1^{h}(f,
\phi))$, where the intersection is taken over all pointed covering
spaces $(f, \phi) \: (\Y,y) \to (\X,x)$.
\end{cor}

\begin{proof}
By Lemma \ref{P6} we have,

\begin{align*}
\bigcap \im(\pi_1^{h}(f, \phi)) &=\bigcap
\omega_x^{-1}(\im(\pi_1(f, \phi))) \\ &= \omega_x^{-1}(\bigcap
\im(\pi_1(f, \phi))) \\ &=\omega_x^{-1}(\{1\})=\Ker \omega_x.
\end{align*}
\end{proof}

\begin{cor}{\label{P7}}
Let $\X$ be an algebraic stack, and let $x$ be a geometric point.
Then the map $\ox \: \phX \to \pX$ has a finite image.
Furthermore, the following conditions are equivalent:

\begin{itemize}
\item[$\mathbf{i})$] The map $\omega_x$ is injective
(respectively, $\omega_x$ is injective and $x$ has a reduced
stabilizer).

\item[$\mathbf{ii})$] There exists a covering  space
$f \: \Y \to \X$ with an unramified (respectively, a schematic)
point $y$ of $\Y$ lying above $x$.
\end{itemize}

\noindent  Furthermore, if $f \: \Y \to \X$ of part
$(\mathbf{ii})$ is Galois, then every point in $\Y$ lying above
$x$ is unramified (respectively, schematic).
\end{cor}

\begin{proof}
That $\ox$ has finite image follows from the fact that $\ox$
factors through the fundamental group of the residue gerbe at $x$,
which is finite (example \ref{P5.9}). The equivalence of
$(\mathbf{i})$ and $(\mathbf{ii})$ follows from Corollary \ref{P5}
and Corollary \ref{P4.9}, plus the fact that $\ox$ has finite
image.  The last statement is obvious.
\end{proof}

\begin{lem}{\label{P7.5}}
If the stabilizer of $\X$ is finite, then the set of all schematic
points of $\X$ is open. This set is the underlying set of the
largest open substack of $\X$ that is an algebraic space.
\end{lem}

\begin{proof}
The first assertion follows from (\cite{EGA4} , 18.12.7) and
Corollary \ref{P4.5}. The second assertion follows from Corollary
\ref{P7} and Lemma \ref{P4.9}.
\end{proof}

\begin{rem}
In fact,  Proposition 18.12.7 of \cite{EGA4} is stated for proper
morphisms, so the above result is still valid if we assume the
stabilizer is only proper.
\end{rem}

\begin{cor}{\label{P7.6}}
Let $\X$ be an algebraic stack with finite stabilizer, and let $f
\: \Y \to \X$ be a representable finite \'{e}tale cover. Let $x$
be a  point in $\X$ so that every points in $\Y$ lying above $x$
is schematic. Then, there exists an open neighborhood $U$ of $x$
so that the inverse image of $U$ in $\Y$ consists entirely of
schematic points (i.e., is an algebraic space).
\end{cor}

\begin{proof}
By Corollary \ref{P4.9} $(\mathbf{ii})$, $\Y$ has finite
stabilizer. Hence, by Corollary \ref{P7.6}, there exists an open
set $V \subseteq \Y$ containing the fiber of $x$ all of whose
points are schematic. Since $f$ is closed, we could find an open
$U \subseteq \X$ whose inverse image is contained in $V$ (say, the
complement of $f(\Y \backslash V)$). The open set $U$ has the
desired property.
\end{proof}

\begin{ex}{\label{E:1}}{\it Counterexample when
stabilizer is not finite.} Take the affine line $\mathbb{A}^1$ and
let the group scheme $\mathbb{A}^1 \coprod \mathbb{A}^1 \backslash
\{0\} \to X$ act (trivially) on it. The resulting quotient stack,
call it $\X$, has only one schematic point. An alternative way of
constructing $\X$ is to mode out the affine line with double
origin by the action of $\mathbb{Z}/2\mathbb{Z}$ that leaves all
the points invariant and swaps the two origins. Hence, $\X$ is
uniformizable and we have $\pi_1(\X) \cong
\mathbb{Z}/2\mathbb{Z}$. By Corollary \ref{P7}, all the maps $\ox$
are isomorphisms, except for the one at the unique schematic point
(``the origin'') for which $\ox$ is the zero map.
\end{ex}

\begin{ex}{\label{E:2}}{\it An algebraic stack with proper
stabilizer where the set of unramified points is not open.} Let
$k$ be a field of characteristic 2. Let $\mathbb{A}^1=\Spec k[z]$
be the affine line, and let $G:=\Spec k[x,y]/{(y^2=yx)}$. Then it
is easy to see that $G$ is a finite flat scheme over
$\mathbb{A}^1$ via the structure map $$\Spec k[x,y]/{(y^2=yx)}
\to     \Spec k[x].$$ A bit less non-trivial is the fact that $G$
is indeed a group scheme over  $\mathbb{A}^1$. This can be seen as
follows. An affine scheme $T$ over $\mathbb{A}^1$ is determined a
pair $(R,r)$ consisting of a $k$-algebra $R$ and an element $r \in
R$; a $T$-point of $G$ over $\mathbb{A}^1$ is determined by an
element $a \in R$ that satisfies the equation $a^2=ar$. This set
is naturally an Abelian group under addition. Therefore, $G$ is
naturally a group scheme over $\mathbb{A}^1$. This group scheme is
a constant group scheme outside the origin, but is non-reduced
over the origin. The quotient of the trivial  action of $G$ on
$\mathbb{A}^1$ is an algebraic stack that has exactly one
unramified point.
\end{ex}

\begin{ex}{\label{E:3}}{\it A Deligne-Mumford stack with finite \'{e}tale
stabilizer where $\ox$ is never injective.} Let $\X$ be the
$\mathbb{Z}/2\mathbb{Z}$-gerbe over $\mathbb{P}^1$ associated to
the nontrivial class in $H^2(\mathbb{P}^1, \mathbb{Z}/2\mathbb{Z})
\cong \mathbb{Z}/2\mathbb{Z}$. For all points $x \in \X$, we have
$\phX \cong \mathbb{Z}/2\mathbb{Z}$. We prove in Section
\ref{S:examples} that $\pX$ is trivial. Therefore, $\ox$ is the
zero map for all points $x$.
\end{ex}

Examples \ref{E:2} and \ref{E:3}  are  simply connected (when the
base field is algebraically closed). For Example \ref{E:2} we can
use an argument as in Example \ref{P5.9}. The proof that Example
\ref{E:3} is simply connected can be found in Section
\ref{S:examples}.
\section{Classification of uniformizable stacks}

The main result of this section is Theorem \ref{T:uniformization}
which gives a necessary and sufficient condition for the
``universal cover'' of an algebraic stack to be an algebraic
space.

\begin{defn}
We say that an algebraic stack $\X$ is {\sl uniformizable} if it
has a covering  space that is an algebraic space.
\end{defn}

By (\cite{LM}, Theorem 6.1), $\X$ is uniformizable if and only if
it is of the form $[X/G]$, where $X$ is an algebraic space and $G$
is a finite group acting on it.

\begin{thm}{\label{T:uniformization}}
Let $\X$ be a Noetherian algebraic stack. Then, $\X$ is
uniformizable if and only if it is Deligne-Mumford and $\ox$ is
injective for any geometric point $x$.
\end{thm}

\begin{proof}
Necessity follows from Corollary \ref{P7}.  Let us prove the
sufficiency. Since the stabilizer of $\X$ is quasi-finite, there
exists a stratification $\emptyset=\X_0 \subset \X_1 \subset
\cdots \subset \X_n=\X$ by close substacks such that the
stabilizer of each $\Z_i=\X_i \backslash \X_{i-1}$ is finite.
Consider $\Z_i$ for some $i$, and let $x$ be a geometric point in
$\Z_i$. Let $f \: \Y \to \X$ be a Galois cover of $\X$ as in
Corollary \ref{P7} (note that we can always replace an arbitrary
cover by a Galois cover that dominates it), and let $V_x \subseteq
\pX$ be the corresponding open normal subgroup (of finite index).
By Corollary \ref{P7.6}, there is an open neighborhood of $V_x$ of
$x$ in $\Z_i$ so that the inverse image of $U_x$ in $\Y$ consists
entirely of schematic points. The $V_x$, for various $x$ in
$\Z_i$, form an open covering of $\Z$, so, by Noetherian
hypothesis, we can find finitely many of $V_x$ that still cover
$\Z_i$. Let $U_i \subseteq \pX$ be the intersection of the
corresponding subgroups. Then, $U_i \subseteq \pX$ is an open
normal subgroup (of finite index). Furthermore, the Galois cover
associated to $U_i$, restricted to $\Z_i$,  consists entirely of
schematic points. Now set $U:=U_1 \cap U_2 \cap \cdots \cap U_n$.
The Galois cover associated to $U$ consists entirely of schematic
points, hence is an algebraic space by Corollary \ref{P4.5}
$(\mathbf{ii})$.
\end{proof}

\begin{ex}
Let $\mathcal{C}$ be a (not necessarily compact) smooth
Deligne-Mumford curve over $\mathbb{C}$ whose moduli space
(necessarily a smooth algebraic curve) is not isomorphic to
$\mathbb{P}^1$, or, if it is, either $\mathcal{C}$ has at least
three corner points (i.e., points at which the stabilizer group
jumps), or it has two corner points with isomorphic stabilizer
groups. In this case, an explicit calculation of the fundamental
groups using van Kampen's theorem shows that all the maps $\ox$
are injective. Therefore, the universal cover of $\mathcal{C}$  is
either $\mathbb{C}$, $\mathbb{P}^1$ or the upper half-plane. For a
complete classification of smooth Deligne-Mumford stack with their
uniformization types and an explicit calculation of their
fundamental groups and maps $\ox$ we refer the reader to
\cite{Unif}.
\end{ex}

\section{Fundamental Group of the Moduli Space}{\label{S:moduli}}

Morally, a {\em moduli space} for an algebraic stack $\X$ is an
algebraic space $\Xm$ that ``best'' approximates the given
algebraic stack. To give a more mathematical sense to this
statement, we would like to have a morphism $\pi \: \X \to \Xm$
and we impose a list conditions on our morphism $f$ to ensure
$\Xm$ is  ``close enough'' to $\X$ in a reasonable sense.
Depending on what conditions we impose on $f$, there is a name
 for the type of moduli space we get. To name a few: categorical
 quotient, uniform categorical quotient,
 topological quotient, geometric quotient moduli, coarse moduli
 space, GC quotient, etc.. Reader can consult
 \cite{Mumford}, \cite{KM}, \cite{Kollar}  for  precise
 definitions. There is a plentiful of results in the literature
 in the form:  {\it For certain a type of algebraic stacks,
 certain type of a moduli space exists.}
 In order not to get bogged down with unnecessary
 hypotheses, in this paper we will not stick to any of these
 standard definitions.  Instead, we list a few
properties that we will need our moduli space to have to make the
proofs work, and we assume we are working in a {\em suitable}
category of algebraic stacks in which such moduli spaces exist
(see Definition \ref{D:suitable} below). It turn out that
basically any reasonable definition of a  moduli space satisfies
our properties. For this reason, in this work, we will use the
term {\em moduli space} without any adjective. We leave it to the
reader to check whether their favorite definition satisfies our
axioms or not.

Vaguely speaking, by a {\em suitable category} $\mathfrak{C}$ of
algebraic stacks we mean a category of algebraic stacks in which a
notion of  {\em moduli space} is specified. Once such a category
is fixed, when we say an algebraic stack $\X$ has a moduli space
we mean $\X \in \mathfrak{C}$. More precisely:

\begin{defn}{\label{D:suitable}}
By a {\sl suitable}\, category $\mathfrak{C}$ of algebraic stacks,
we mean a subcategory of the category of algebraic stacks that
contains the category of algebraic spaces, and in which for every
algebraic stack $\X$ we have chosen an algebraic space $\Xm$
(called the {\sl moduli space} of $\X$), and a morphism $\pi_{\X}
\: \X \to \Xm$ (called the {\sl moduli map}), that satisfy the
following axioms:

\begin{itemize}
\item[$\mathbf{M1.}$] {\it Functoriality.} If $f \: \Y \to \X$ is
(2-isomorphism class of) a morphism of stacks, then we have a
morphism of moduli spaces $f_{mod} \: \Ym \to \Xm$ which makes the
following diagram commute:

$$\xymatrix{ \Y \ar[r]^f \ar[d]_{\pi_{\Y}}  &   \X
\ar[d]^{\pi_{\X}} \\
           \Ym \ar[r]_{f_{mod}}           &   \Xm}$$
   with the usual functoriality properties (i.e., $\X \mapsto \Xm$ is
   a functor and $\pi$ is a natural transformation
   between the identity functor and the moduli functor).

\item[$\mathbf{M2.}$]{\it Geometric points.} For any algebraically closed
field $k$, $\pi_{\X}$ induces a bijection
 between $k$-points (up to 2-isomorphism) of $\X$ and $k$-points of $\Xm$.

\item[$\mathbf{M3.}$]{\it Covering spaces.} If $\X$ is in $\mathfrak{C}$,  then
so is every  covering  space  of $\X$.

\item[$\mathbf{M4.}$] {\it Invariance under finite \'{e}tale base change.}
If we have a cartesian diagram

$$\xymatrix{ \Y \ar[r]^f \ar[d]_{\pi}  &   \X \ar[d]^{\pi_{\X}} \\
           Y \ar[r]_g           &   \Xm,}  $$
in which $g$ is finite \'{e}tale, then $Y$ is the moduli space of
$\Y$ with  $\pi$ its   moduli map. Furthermore, $g$ coincides with
the induced map of moduli spaces as in $\mathbf{M1}$.

\item[$\mathbf{M5.}$] {\it Free quotient of a moduli map is a moduli map.}
Let $G$ be a finite group acting freely (see Definition
\ref{D:free} below) on $\Y$, and assume $\Y/G$ is in
$\mathfrak{C}$. Then in the following (necessarily cartesian)
diagram, $\pi_{\Y}/G$ is a moduli map, and the quotient map $q$ is
the same as the induced map in $\mathbf{M1}$.

$$\xymatrix{ \Y \ar[d]_{\pi_{\Y}} \ar[r]^f & \Y/G
\ar[d]^{\pi_{\Y}/G} \\
             \Ym \ar[r]_g                & \Ym/G. }$$

\end{itemize}
\end{defn}

\noindent {\bf Convention.} Throughout this paper, whenever we
talk about the moduli space of an algebraic stack $\X$, it is
assumed that $\X$ belongs to a certain fixed suitable category as
defined above.

\begin{defn}{\label{D:free}}
Let $G$ be a finite group acting (via
morphisms-up-to-2-isomorphism) on  an algebraic stack $\X$, and
let $x$ be a geometric point. We define the {stabilizer group} of
$x$ to be the set of all $g \in G$ such that $g(x)$ is
2-isomorphic to $x$. We say that the action is {\it free}, if the
stabilizer group of every  geometric point $x$ is trivial.
\end{defn}

If $\X$ has a moduli space that satisfies $\mathbf{M1}$ and
$\mathbf{M2}$, then stabilizer group of a point is equal to the
usual stabilizer group of the corresponding point in $\Xm$ under
the induced action of $G$. Therefore, in this case, the action
being free means that the induced action of $G$ on $\Xm$ is free.

A caveat about the above definition is in order. Unlike the case
of algebraic spaces, the Galois group of a (representable) Galois
cover $\Y \to \X$ may {\em not} act freely on $\Y$ (although it
acts freely on the fibers $F_x(\Y)$).

\begin{lem}{\label{L:free}}
Let $\X$ be an algebraic stack, and let $f \: \Y \to \X$ be a
Galois cover with Galois group $G$. Let $y$ be a geometric point
of $\Y$ and $x=f(y)$ its image in $\X$. Then we have a short exact
sequence

$$1 \to \phY \to \phX \to G_y \to 1.$$
\end{lem}

\begin{proof}

Exactness on the left follows from Lemma \ref{D1}. Since $\Y \to
\X$ is Galois, the action of an element $\gamma \in \phX$ (more
precisely, its image in $\pX$ via $\ox$) on the fiber extends to
an element $g_{\gamma} \in G$. Recall from Section \ref{S:Galois}
that this action takes $(y, id) \in F_x(\Y)$ to $(y, \gamma^{-1})$
(here by $id$ we mean the constant transformation from $x$ to
itself). Therefore, $g_{\gamma}$ indeed belongs to $G_y$. This
defines a natural group homomorphism $\phX \to G_y$. Proof of
surjectivity of this map is a easy, but a bit intricate. We need
to get down to precise definitions. Recall from the definition of
the Galois category $\mathbf{C}_{\X}$ (Section \ref{S:Galois} that
an element in $g \in G$ can be represented by a pair $(a,\Phi)$,
where $\Phi \: f \Ra f \circ a$. This element belongs to $G_y$ if
and only if   there is a 2-isomorphism $\alpha \: y \rsa
a(y)$. The effect of $(a,\Phi)$ takes $(y, id) \in F_x(\Y)$ to
$(a(y), \phi) \sim (y,\phi f(\alpha^{-1})) \in F_x(\Y)$, where
$\phi \: x \rsa f(a(y))$ is the transformation induced by
$\Phi$. We claim that $\gamma=(\phi f(\alpha^{-1}))^{-1} \in \phX$
maps to $g$, i.e., we have $g=g_{\gamma}$. But we just verified
that these two cover automorphisms have the same effect on $(y,
id) \in F_x(\Y)$, so they must be equal.  This proves
surjectivity. Exactness in the middle follows from a similar
`hidden path chasing' argument.
\end{proof}

\begin{prop}
Let $\X$ be an algebraic stack, and let $f \: \Y \to \X$ be a
Galois cover with Galois group $G$. Then the following are
equivalent:

\begin{itemize}
\item[$\mathbf{i}$)] The action of $G$ on $\Y$ is free
(see Definition \ref{D:free}).

\item[$\mathbf{ii}$)] For any geometric point $y$, the induced
map $\phY \to \pi_1^h(\X,f(y))$ is an isomorphism. Equivalently,
the induced map of of groupoids $\PhY \to \PhX$ is fully faithful.

\item[$\mathbf{iii}$)] Let $X \to \X$ be a chart for $\X$ and
$Y \to \Y$ the pull back chart for $\Y$. Then the morphism $S_Y
\to S_X \x_X Y$ of Lemma \ref{P4.8} is an isomorphism.
Equivalently, the pull back of the stabilizer group stack of $\X$
is the stabilizer group stack of $\Y$ (see Section \ref{S:Stab}).
\end{itemize}
\end{prop}

\begin{proof} The equivalence of  $(\mathbf{i})$ and  $(\mathbf{ii})$
is immediate from Lemma \ref{L:free}. To prove the equivalence of
$(\mathbf{ii})$ and  $(\mathbf{iii})$, note that, by Corollary
\ref{P4.9}, the map $S_Y \to S_X \x_X Y$ is an open and closed
immersion. The equivalence follows immediately from Proposition
\ref{P4}.
\end{proof}

In fact, the equivalence of ($\mathbf{ii}$) and ($\mathbf{iii}$)
is true for any covering map $f \: \Y \to \X$ (not necessarily
Galois).

\begin{defn}{\label{D:FPR}}
When the two equivalent conditions alluded to in the previous
paragraph are satisfied for a covering map $f$, we say that $f$ is
a {\em fixed point reflecting} (or FPR) morphism.
\end{defn}

The following lemma provides a source of FPR morphisms.

\begin{lem}{\label{G7}}

Let $\pi \: \X \to A$ be an arbitrary map to an algebraic space,
and let $f \: B \to A$ be an arbitrary map of algebraic spaces.
Then, the induced map $\Y :=\X \x_A B \to \X$ is FPR.

\end{lem}

\begin{proof}

Let $p \: X \to \X$ be a chart for $\X$ and $Y:=X\x_{\X}\Y=X\x_AB$
the corresponding pull back chart for $\Y$. It follows from the
definition of the stabilizer group that the diagram

$$\xymatrix{ S_Y \ar[r] \ar[d] & S_X \ar[d] \\
             Y \ar[r] & X }$$
 is cartesian. The lemma follows now from Proposition \ref{P4}.
 \end{proof}

\begin{prop}{\label{G8}}

Let $\X$ and $\Y$ be algebraic stacks with $\X$ connected. Let $f
\: \Y \to \X$ be an FPR covering space. Then, the induced map on
the moduli spaces $f_{mod} \: \Ym \to \Xm$ is finite \'{e}tale.
Furthermore, the following diagram is cartesian:

$$\xymatrix{ \Y \ar[r]^{f} \ar[d] & \X \ar[d] \\
             \Ym \ar[r]_{f_{mod}} & \Xm }$$
\end{prop}

\begin{proof}
We will only need the result when $f$ is Galois, in which case the
claim is immediate from $\mathbf{M5}$. (Hint: for the general case
take a Galois cover $\Y' \to \X$ that factors through $\Y \to \X$,
and use the fact that $\Y' \to \Y$ is also Galois.)
\end{proof}

\begin{prop}{\label{G9}}
Let $\X$ be an algebraic stack. Then,  we have a one-to-one
correspondence $$\xymatrix@C=60pt{ \left\{ \txt\footnotesize{
Fixed point reflecting covering spaces of $\X$.}\right\} \ar
@/^/[r]^(0.57){\txt\tiny{\bf{moduli}}}  &
\left\{\txt\footnotesize{Covering spaces of $\Xm$. }\right\}  \ar
@/^/ [l]^(0.42){\txt\tiny{\bf{base extension}\\ {\bf via}
$\pi_{mod}$}}}.$$ \noindent Similar statement is true for
connected covering spaces. Finally, the lemma remains valid if we
assume everything  pointed.
\end{prop}

\begin{proof}
The assertion is an immediate consequence of Lemma \ref{G7} and
Lemma \ref{G8}, plus the axiom $\mathbf{M4}$. Also note that in a
diagram as in Lemma \ref{G8}, $\Y$ is connected if and only if
$\Ym$ is so. The pointed version of the lemma is true because a
moduli map induces a bijection on the geometric points
($\mathbf{M2}$).
\end{proof}

\begin{lem}{\label{G10}}

Let $G$ be a profinite group and $H \subseteq G$ a normal
subgroup. Then, the intersection of all open normal subgroups of
$G$ that contain $H$ is equal to $\bar{H}$.
\end{lem}

\begin{proof}

The intersection obviously contains $\bar{H}$. To prove the
equality, we may assume $H$ is closed. Take an element $g \in G
\backslash H$. We want to find an open normal subgroup of $G$ that
contains $H$ but not $g$. Since $H$ is closed, and since open
normal subgroups of $G$ form a fundamental system of neighborhoods
at the origin, we can find an open normal subgroup $U \subset G$
such that $gU \cap H$ is empty. The group $HU$ is easily seen to
have the desired property.
\end{proof}

\begin{lem}{\label{G10.5}}

Let $\mu \: G_1 \to G_2$ be a continuous homomorphism of profinite
groups. Let $N \subseteq G_1$ be a closed normal subgroup. Assume
$\mu^{-1}$ induces a one-to-one correspondence between open
subgroups of $G_2$ and those open subgroups of $G_1$ that contain
$N$. Then,  $\mu$ induces an isomorphism $\bar{\mu} \: G_1/N
\risom G_2$.
\end{lem}

\begin{proof}

The intersection of all open normal subgroups of $G_2$ is $\{ 1
\}$; so, by the above correspondence, the intersection of all open
normal subgroups of $G_1$ that contain $N$ is equal to
$\mu^{-1}(\{1\})$. So Lemma \ref{G10} implies that
$\mu^{-1}(\{1\})=N$. Hence, $\bar{\mu}$ is injective. To prove the
surjectivity, consider the closed subgroup $H:=\mu (G_1)$ of
$G_2$. Assume $H \neq G_2$. Then, by Lemma \ref{G10}, there exists
a proper open subgroup $U$ of $G_2$ that contains $H$. But $U$ and
$G_2$ are now two different open subgroups of $G_2$ whose inverse
image under $\mu$ is equal to $G_1$; a contradiction.
\end{proof}

Let $\X$ be a connected algebraic stack that has a moduli space
(Definition \ref{D:suitable}). Fix a geometric point $x$ for $\X$.
For any other geometric point $x'$, and any choice of a path from
$x'$ to $x$, we may identify the image of $\omega_{x'}$ in
$\pi_1(\X,x')$ with a subgroup of $\pX$. Let $N \subseteq \pX$
denote the closure of the subgroup generated by all these
subgroups (for various $x'$ and various paths connecting $x'$ to
$x$). Then, $N$ is a normal subgroup of $\pX$ that maps to zero
via the homomorphism $\pX \to \pi_1(\Xm,x)$.

\begin{thm}{\label{G5}}
 The natural map $\pX/N \to \pXm$ is an isomorphism.
\end{thm}

\begin{proof}
We may assume that $\X$ is connected.  We want to show that the
map $\mu \: \pX \to \pXm$ induces an isomorphism $$\bar{\mu} \:
\pX/N \risom \pXm .$$

By Lemma \ref{G10.5}, it is enough to show that $\mu^{-1}$ induces
a one to one correspondence between open subgroups of $\pXm$ and
those open subgroups $U$ of $\pX$ that contain $N$. If we unravel
the definition of $N$, we see that for $U$ to contain  $N$ is
equivalent to the following: For any geometric point $x'$ and any
path $\gamma$ from  $x$ to $x'$, the image of $U$ under the
isomorphism $\pX \risom \pi_1(\X,x')$ induced by $\gamma$ contains
$\pi_1^{h}(\X,x')$. By Lemma \ref{P6}, this is equivalent to $\X_U
\to \X$, the corresponding (pointed) covering  space, being FPR.
This shows  that there is a one-to-one correspondence between open
subgroups of $\pX$ that contain $N$ and (pointed) connected FPR
covering  spaces of $\X$. By Lemma \ref{G9}, the latter set is,
via base change along $\pi_{mod} \: \X \to \Xm$ , in one-to-one
correspondence with (pointed) connected covering spaces of $(\Xm,
x)$. Finally, this last set is in one-to-one correspondence with
open subgroups of $\pXm$. Proof is now complete.
\end{proof}

\section{Application}{\label{S:application}}

In this section we apply the formula of Theorem \ref{G5} to
calculate the fundamental group of the G.I.T. quotient of a
groupoid actions. We will make use of the notion of a {\em
fibration} as sketched in  the appendix.  All schemes (algebraic
spaces, stacks) and all morphisms are assumed to be quasi-compact.
{\bf Notation:} When $R \st{} X$ is a groupoid, the notation $X/R$
refers to the quotient as an {\em algebraic space} (if it exists),
as opposed to $[X/R]$, which denotes the quotient as a stack. More
precisely, $X/R$ denotes the {\em moduli space} of $[X/R]$ (in the
sense of Section \ref{S:moduli}), if it exists.

The morphisms that we call {\em fibrations} have the following
properties (see the appendix for details):

{\small
\begin{itemize}

\item[$\triangleright$] Let $(f,\phi) \: (\Y,y) \to (\X,x)$ be a pointed
fibration. Then we have an exact sequence \\ $\pi_1(\Y_x,y) \to
\pY \to \pX  \to  \pi_0(\Y_x,y) \to  \pi_0(\Y,y) \to \pi_0(\X,x)
\to \{*\}.$

\item[$\triangleright$] Composition of fibrations is again a fibration.
An arbitrary base extension of a fibration is a fibration.

\item[$\triangleright$] Being a fibration is local (on the target) in the fppf topology.

\item[$\triangleright$] A morphism $f$ is a fibration if and only if $f_{red}$ is a fibration.
\end{itemize}
}

We just mention a few examples of fibrations:

\begin{itemize}

\item[$\mathbf{i}$)] Proper flat morphisms that have geometrically reduced fibers (in particular,
finite \'{e}tale morphisms).

\item[$\mathbf{ii}$)] Any morphism whose target is spectrum of a field (or spectrum of an Artin ring),
and any base extension of such a morphism.

\item[$\mathbf{iii}$)] The structure map of a smooth group scheme whose geometric fibers have a fixed
number of connected components.

\item[$\mathbf{iv}$)] In a geometric context (e.g., when everything is over $\mathbb{C}$), a morphism
$f$ that is a (quasi-)fibration in the topological sense is a
fibration in our sense.
\end{itemize}

Let $X$ be a connected scheme (or algebraic space) and let $R
\st{s,t} X$ be a groupoid which is flat and of finite presentation
and whose diagonal map $R \to X \x X$ is separated (so the
quotient stack becomes algebraic). Assume further that $s$ (hence
$t$) is a fibration. Fix a geometric point $x \: \Spec k \to X$,
and let $S_x$ be the stabilizer group at $x$. Let $R_x$ be the
fiber of $s \: R \to X$ over $x$. There is a natural base point
for $R_x$, namely the image of $x$ under the identity section $X
\to R$. It follow from the definition of a groupoid, that the
group of $k$-points of the stabilizer group scheme at $x$ (which,
by abuse of notation, we call $S_x$) acts on the  set of
$k$-points of $R_x$, and, a fortiori,
 on the  set $\mathbf{R}_x$ of connected components of  $R_x$. Let $x' \: \Spec k' \to X$
be another geometric point. Any choice of a ``path'' from $x'$ to
$x$ induces an action of $S_{x'}$ on $\mathbf{R}_x$ (because,
since $R \to X$ is a fibration, $\mathbf{R}_x$ form a
$\Pi_1(X)$-set). Therefore, we have an action of $S_{x'}$
 on $\mathbf{R}_x$ that is well-defined up to conjugation by
the action  of $\pi_1(X,x)$ on $\mathbf{R}_x$. Let $M_x$ be the
largest quotient of $\mathbf{R}_x$ on which every $S_x$ acts
trivially (in other words, the collection of  $M_x$ forms a
$\Pi_1(X)$-set that is the largest quotient of the  $\Pi_1(X)$-set
\{$\mathbf{R}_x$\} on which the action of all $S_x$'s is trivial).
There is an alternative way of constructing $M_x$: Let $\X=[X/R]$,
and consider the fibration $p \:X \to \X$ (See Theorem
\ref{T:last}). Let $$\pi_1(X,x) \to \pi_1(\X,p(x)) \to
\pi_0(X_{p(x)},x) \to \{*\}$$ be the corresponding fiber homotopy
exact sequence. The (pointed) set  $\pi_0(X_{p(x)},x)$ is
canonically isomorphic to the (pointed) set of connected
components of $R_x$. Let $N \subseteq \pi_1(\X,p(x))$ be as in
Theorem \ref{G5}. If we now kill the action of $N$ on
$\pi_0(X_{p(x)},x)$, what  we obtain is canonically isomorphic (as
a (pointed) set with a $\pi_1(X,x)$ action) to $M_x$. Theorem
\ref{G5} implies the following

\begin{thm}{\label{T:quot1}}
Let $R \st{} X$ be a flat of finite presentation groupoid such
that the diagonal $R \to X \x X$ is separated ($R$ and $X$
algebraic spaces) and let $R/X$ be a quotient (i.e., a moduli
space in the sense of Section \ref{S:moduli} for the algebraic
stack [X/R]), if it exists. Let $q \: X \to X/R$ be the quotient
map. Then we have a canonical exact sequence $$\pi_1(X,x) \to
\pi_1(X/R,q(x)) \to M_x \to \{*\}.$$
\end{thm}

Now, let $f \: X \to S$ be a connected algebraic space over a
connected base $S$, and let $G$ be a flat, separated and of finite
presentation group space over $S$ acting on $X$. Assume further
that $G \to S$ is a fibration. Suppose a quotient $X/G$ exists (in
the sense of the above theorem). Fix a geometric point $x \: \Spec
k \to X$. Let $s=f(x)$. Let $\mathbf{G}$ denote the set of
connected components of $G_s$. Define $\mathbf{H}_x \subset
\mathbf{G}$ to be the set of all connected components of $G_s$
that contain at least one $k$-point that leaves $x$ fixed. For any
other choice of a geometric point $x'$
 of $X$, and for any choice of a path from $x'$ to $x$, we can
 identify $\mathbf{H}_{x'}$ with a subgroup of
$\mathbf{G}$. Let $\mathbf{I} \subset \mathbf{G}$ be the subgroup
generated by all these groups. In this situation, Theorem
\ref{T:quot1} can be translated as follows

\begin{thm}{\label{T:quot2}}
Let $G$ be a group space acting on a connected algebraic spaces
$X$, relative to a connected base $S$. Suppose $G \to S$ is a
fibration that is flat, separated and of finite presentation.
Assume a quotient $X/S$ exists (in the sense of Theorem
\ref{T:quot1}. Then, we have an exact sequence $$\pi_1(X) \to
\pi_1(X/G) \to \mathbf{G/I} \to \{*\}.$$ Here, $G/I$ is viewed as
a pointed set.
\end{thm}

In the above situation,  if $X \to S$ has connected set-theoretic
(as opposed to geometric) fibers, then $\mathbf{I} \subset
\mathbf{G}$ could also be defined as follows: Let $\mathbf{H}
\subset \mathbf{G}$ be the set of all connected components of
$G_s$ that contain at least one $k$-point whose action on $X_s$
has a fixed point. For any other choice of a geometric point $s'$
that lies in the image of $X$, and for any choice of a path from
$s'$ to $s$ coming from $X$, we can identify $\mathbf{H'}$ (the
counterpart of $\mathbf{H}$ at $s'$) with a subgroup of
$\mathbf{G}$. Then $\mathbf{I} \subset \mathbf{G}$ is the subgroup
generated by all these groups. Also, the same is description for
$\mathbf{I}$ is valid when the action of $\pi_1(X)$ (via that of
$\pi_1(S)$) on the set of connected components of the geometric
fibers of $G \to S$ is trivial. In this case, $\mathbf{I} \subset
\mathbf{G}$   will indeed be a normal subgroup and the exact
sequence of Theorem \ref{T:quot2} will become an exact sequence of
groups.  The following corollaries are instances where this
phenomenon occurs.

\begin{cor}{\label{C:quot3}}
Let $X$ be a connected algebraic space, and let $G$ be a finite
group acting on $X$ (relative to a certain base $S$). Let $X/G$ be
the quotient for this action. Let $I \subset G$ be the subgroup
generated by all elements that have fixed points. Then, $I$ is
normal in $G$, and we have the following exact sequence of groups:
$$ \pi_1(X) \to \pi_1(X/G) \to G/I \to 1.$$
\end{cor}

\begin{cor}{\label{C:quot4}}
Let $X$ be a connected algebraic space over an algebraically
closed field $k$, and let $G$ be an algebraic group (not
necessarily reduced) over $k$ acting on $X$. Assume a quotient
$X/G$ exists (in the sense of Theorem \ref{T:quot1}). Let
$\mathbf{G}$ be the group of connected components of $G$, and let
$\mathbf{I} \subset \mathbf{G}$ be the set of all components that
contain at least one element that leaves some point of $X$ fixed.
Then, $\mathbf{I}$ is a normal subgroup of $\mathbf{G}$, and we
have the following exact sequence of groups: $$ \pi_1(X) \to
\pi_1(X/G) \to \mathbf{G/I} \to 1.$$
\end{cor}

\begin{cor}{\label{C:quot5}}
Let $X$ and $G$ be as in Theorem \ref{T:quot2}. Assume further
that $X$ is simply connected. Then, we have an isomorphism $
\pi_1(X/G) \risom \mathbf{G/I}$ .
\end{cor}

\subsection{Description of the kernel of $\pi_1(X) \to \pi_1(X/G)$}

 If we go through the steps that resulted in the proof of Theorem
 \ref{T:quot1}, we will see that the kernel of $\pi_1(X) \to
 \pi_1(X/G)$ is, vaguely speaking, the smallest subgroup of
 $\pi_1(X)$ containing all the elements that ``obviously'' map to
 zero. To illustrate this idea, let us consider the typical
 situation  where $G$ is a discrete group acting on a topological
 space.\footnote{The results of this paper can be naturally
 generalized to topological/combinatorial stacks, in which case
 instead of the algebraic fundamental group we work with the usual
 fundamental group. The topological approach is the subject of a
 forthcoming paper.} Let $X/G$ be the usual topological quotient,
 and let $q \: X \to X/G$ denote the quotient map. Pick an element
 $g \in G$ and a fixed point $y \in X$ for $g$. Let $x \in Y$ be
 an arbitrary point and $\gamma$ a path joining $x$ to $y$. Then
 $\gamma g(\gamma^{-1})$, the juxtaposition of $\gamma$ and
 $g(\gamma^{-1})$, is a path from $x$ to $g(x)$, whose image in
 $X/G$ is a trivial loop at $q(x)$. Let us call a path that is of
 the form  $\gamma g(\gamma^{-1})$, or more generally any path
 that can be obtained by composing a finite number of such paths,
 a {\em doomed} path. In particular, we can talk about doomed
 loops around $x$, namely, doomed paths which  start from $x$ and
 end at $x$. Doomed loops form a normal subgroup of $\pi_1(X)$,
 and they obviously map to zero in $\pi_1(X/G)$. It follows from
 the proof of Theorem \ref{T:quot1} that, in fact, the group of
 doomed loops is equal to the kernel of $\pi_1(X) \to \pi_1(X/G)$.

 When $G$ is no longer discrete, there are some more loops that
 ``obviously'' map to zero in $\pi_1(X/G)$; namely, the ones
 coming from the loops in $G$. More precisely, let $x$ be a point
 in $X$. Then, we have a map from $G$ to $X$ that sends $g$ to
 $g(x)$. The image of the fundamental group of (the connected
 component of the identity of) $G$ in $\pi_1(X,x)$ lies in the
 kernel of $\pi_1(X) \to \pi_1(X/G)$. Again, it follows from the
 proof of Theorem \ref{T:quot1} that this group, together with the
 group of doomed loops, generate the kernel of $\pi_1(X) \to
 \pi_1(X/G)$.

\section{More examples}{\label{S:examples}}

\begin{ex}{\em An example from number theory.}
This example is  to illustrate the relationship between the notion
of ramification on an algebraic stack to the classical notion of
ramification in algebraic number theory. Let $K \subset L$ be a
Galois extension of number fields with Galois group $G$. Let
$\mathcal{O}_K \subset \mathcal{O}_L$ be the corresponding rings
of integers. Let $\X=[\Spec \mathcal{O}_L/G]$. The moduli space of
$\X$ is $\Spec \mathcal{O}_L/G=\Spec \mathcal{O}_K$. The ramified
points of $\X$ corresponds to  primes of $\mathcal{O}_K$ that are
ramified in $\mathcal{O}_L$, and the hidden fundamental group at
each ramified point is isomorphic to the inertia group at the
corresponding prime. Therefore, the group $I \subseteq G$
introduced in Corollary \ref{C:quot3} is equal to the subgroup
generated by all the inertia groups for various primes in
$\mathcal{O}_K$. The exact sequence of Corollary \ref{C:quot3} now
takes the form $$\pi_1(\Spec \mathcal{O}_L) \to \pi_1(\Spec
\mathcal{O}_K) \to G/I \to 1,$$ \noindent which is the same as the
exact sequence $$\Gal(\mathcal{O}_L^{ur}/\mathcal{O}_L) \to
\Gal(\mathcal{O}_K^{ur}/\mathcal{O}_K) \to G/I \to 1,$$ \noindent
whose validity can be checked readily by easy number theoretic
arguments. Here, $\mathcal{O}_L^{ur}$ stands for the maximal
unramified extension.
\end{ex}

\begin{ex}{\label{E:proj}}{\sl Weighted projective spaces are
simply connected.} Let $k$ be an algebraically closed field. Let
$n_1,n_2,\cdots,n_d$ be a sequence of integers. Consider the
action of $k^*$ on $\mathbb{A}^{d}\backslash\{0\}$ which, for any
$\zeta \in k^*$, sends $(x_1,x_2,\cdots,x_d) \in
\mathbb{A}^{k}\backslash\{0\}$ to
$(\zeta^{n_1}x_1,\zeta^{n_2}x_2,\cdots,\zeta^{n_d}x_d)$. Let
$\mathbb{P}(n_1,n_2,\cdots,n_d):=\mathbb{A}^{d}\backslash\{0\}/{k^*}$
be the quotient of this action. It is an algebraic variety   which
is called the {\sl weighted projective space} of weight
$(n_1,n_2,\cdots,n_d)$. It is the moduli space of
$\mathcal{P}(n_1,n_2,\cdots,n_d):=[\mathbb{A}^{d}\backslash\{0\}/{k^*}]$.
The quotient map $\mathbb{A}^{d}\backslash\{0\} \to
\mathcal{P}(n_1,n_2,\cdots,n_d)$ is a fibration. An easy fiber
homotopy exacts sequence argument shows that
$\mathcal{P}(n_1,n_2,\cdots,n_d)$ is simply connected. Therefore,
$\mathbb{P}(n_1,n_2,\cdots,n_d)$ is also simply connected by
Theorem \ref{G5}. \end{ex}

\begin{ex}{$\bar{\mathcal{M}}_{1,1}$.}
It is well-known that the compactified moduli stack
$\bar{\mathcal{M}}_{1,1}$ of elliptic curves (over, say,
$\mathbb{C}$) is isomorphic to $\mathbb{P}(4,6)$, which is simply
connected. Its moduli space is the $j$-line, which is simply
connected, as anticipated by Theorem \ref{G5}.
\end{ex}

\begin{ex}{$\mathcal{M}_{1,1}$.}
The non-compact moduli stack of elliptic curves is the quotient of
the upper half-plane by the action of $SL_2(\mathbb{Z})$.
Therefore, its (topological) fundamental group is isomorphic to
$SL_2(\mathbb{Z})$, and its algebraic fundamental group is
isomorphic to $\widehat{SL_2(\mathbb{Z})}$, the profinite
completion of $SL_2(\mathbb{Z})$. The standard generators $S$ and
$T$ together with $-I$ generate $SL_2(\mathbb{Z})$, and they all
have fixed points. Therefore, by Corollary \ref{C:quot5}, the
quotient of the upper half-plane by the action of
$SL_2(\mathbb{Z})$ (i.e., the moduli space of $\mathcal{M}_{1,1}$)
should be simply connected, which is obviously the case, since it
is isomorphic to $\mathbb{C}$. (To be more precise, in order to be
able to  make use of Corollary \ref{C:quot5}, we should first take
a normal subgroup of $SL_2(\mathbb{Z})$ that acts
fixed-point-freely on the upper half-plane, and then apply
Corollary \ref{C:quot5} to the action of the cokernel of this
group on the quotient space, which is now an algebraic curve.
Otherwise, we need to use the topological counterpart of Corollary
\ref{C:quot5}, which is indeed true.)
\end{ex}

Next we show that the stack of Example \ref{E:3} is
simply-connected. For this, we use the following result from
\cite{Noohi} (also see Corollary \ref{C:seq}).

\begin{thm}{\label{G11}}
Let $\X$ be a connected Deligne-Mumford gerbe with finite
(\'{e}tale) stabilizer. Then, for every geometric point $x$ of
$\X$, the sequence $$\pi_1^{h}(\X,x) \to \pX \to \pXm \to 1  $$
\noindent is exact. If this sequence is exact on the left for some
$x$, then $\X$ is uniformizable.
\end{thm}

\begin{ex}({\em Example \ref{E:3} revisited.})
We have $\phX \cong \mathbb{Z}/2\mathbb{Z}$ for all points $x \in
\X$. Since  $\mathbb{P}^1$ is simply connected, $\X$ does not
admit any  non-trivial finite \'{e}tale cover by a scheme, nor
even by an algebraic space (\cite{Kn}, Corollary 6.16), because,
otherwise, composing it with $\pim$ would give us a non trivial
finite \'{e}tale cover of $\mathbb{P}^1$ (note that the moduli map
$\X \to \mathbb{P}^1$ is finite \'{e}tale). Therefore, $\X$ is not
uniformizable. So by Theorem \ref{G11}, all the maps $\omega_x \:
\phX \to \pX$ are zero maps and $\pX$ is trivial. In other words,
$\X$ does not admit any non-trivial representable finite \'{e}tale
cover.
\end{ex}

The next pathological example is a good exercise to work out the
ideas presented in this paper. It was suggested to me by J. de
Jong.
\begin{ex}{\label{E:pathological}}({\em A non simply-connected,
non-uniformizable stack.})
 Let $N \subset \mathbb{P}^3_{\mathbb{C}}$
be the hypersurface defined by the homogeneous equation
$X^4+Y^4+Z^4+W^4=0$. Let $\zeta$ be a primitive fourth root of
unity. The action of $\mathbb{Z}/4\mathbb{Z}$ on
$\mathbb{P}^3_{\mathbb{C}}$ given by $$(x:y:z:w) \mapsto (x:\zeta
y:\zeta^{2}z:\zeta^{3}w)$$ induces a fixed point free action on
$N$. Let $M$ be the quotient of this action. The fundamental group
of the scheme $M$ is isomorphic to $\mathbb{Z}/4\mathbb{Z}$. Let
$I=\mathbb{Z}/3\mathbb{Z}$. From the spectral sequence $$
H^q(\pi_1(M),H^p(N,I)) \Rightarrow H^{p+q}(M,I)$$ \noindent we
obtain the following exact sequence: $$H^2(\pi_1(M),I)
\hookrightarrow H^2_{et}(M,I) \to H^2_{et}(N,I)^{\pi_1(M)}\to
H^3(\pi_1(M),I).$$ \noindent Note that the action of $\pi_1(M)$ on
$I$ is trivial. So, we have
$H^2_{et}(N,I)^{\pi_1(M)}=H^2_{et}(N,I)$. On the other hand,
$H^3(\pi_1(M),I)=H^3(\mathbb{Z}/4\mathbb{Z},\mathbb{Z}/3\mathbb{Z})=0$.
Therefore, the above exact sequence looks as follows: $$0 \to
H^2(\pi_1(M),I) \to  H^2_{et}(M,I) \to H^2_{et}(N,I) \to  0.$$
\noindent Since $N$ is projective,  $H^2_{et}(N,I)$ is
non-trivial. In particular, $H^2_{et}(M,I)$ is strictly bigger
than $H^2(\pi_1(M),I)$. It is shown in (\cite{Noohi}, Section 6)
that an element in $H^2_{et}(M,I)$ that is not in the image of
$H^2(\pi_1(M),I)$ corresponds to a gerbe $\X$ over $M$ that is
{\em non-uniformizable}. The automorphism group of this gerbe is
$\mathbb{Z}/3\mathbb{Z}$. By Theorem \ref{G11}, we have the
following exact sequence that is {\it not} left exact:
$$\mathbb{Z}/3\mathbb{Z} \to \pi_1(\X) \to \mathbb{Z}/4\mathbb{Z}
\to 0.$$ \noindent Therefore, $\pi_1^h(\X) \to \pi_1(\X)$ is the
zero map, and $\pi_1(\X) \risom \pi_1(\Xm)=\mathbb{Z}/4\mathbb{Z}$
by Theorem \ref{G5}.
\end{ex}

\pagebreak

{\Large\part{Further results}}{\label{C:further}}

\vspace{0.5in}

 Part two of this paper is in a way a technical companion to Part
 one. The first main result in this part (Theorem \ref{T:main}),
 which is a generalization of Theorem \ref{T:uniformization} of
 Part one, tells us to what extent one can resolve stackyness
 after by passing to covering spaces. The second main result is a
 finiteness theorem about the kernel of the homomorphism $\pX \to
 \pXm$ (Theorem \ref{T:finite}). The main tool in proving these
 result is the Stratification Theorem \ref{T:strat}, which we
 believe is useful in its own right.

Let us first reformulate the results of Section \ref{S:basic} in a
form that is more suitable for our next applications:

\begin{prop}{\label{P:a}}
Let $(\X,x)$ be a connected algebraic stack. Then, the
(isomorphism classes of) pointed covering spaces $(\Y,y)$ of
$(\X,x)$ that are FPR at $y$, are in one-to-one correspondence
with open-closed subgroups of $\pX$ that contain the image of $\ox
\: \phX \to \pX$. The (isomorphism classes of) pointed covering
spaces $(\Y,y)$ of $(\X,x)$ that are FPR at every $y'$ in the
fiber of $x$, are in one-to-one correspondence with  open-closed
subgroups of $\pX$ that contain the  normal closure of the image
of $\ox \: \phX \to \pX$.
\end{prop}

If $x'$ is another point in $\X$ (possibly identical to $x$), we
can give a similar description of pointed covers of $(\X,x)$ that
are FPR for all point above $x'$. To do so, we take the normal
closure of the image of $\omega_{x'}$ in $\pi_1(\X,x')$, and
transfer it into a subgroup  of $\pX$ via a path joining $x'$ and
$x$. This subgroup is well defined and can be identified with
$\Phi_{x'}(x)$, where $\Phi_{x'}$ is the normal subgroupoid  of
$\PX$ generated by the image of $\omega_{x'}$ (See Section
\ref{S:groupoid} for definitions and notations). The pointed
covering spaces $(\Y,y)$ of $(\X,x)$ that are FPR at every point
above $y'$ are in one-to-one correspondence with open-closed
subgroups of $\pX$ that contain $\Phi_{x'}(x)$. Indeed, the same
thing is true if we replace $x'$ with a set $S$ of points of $\X$:

\begin{prop}{\label{P:b}}
Let $(\X,x)$ be a connected pointed stack. Let $S$ be a set of
geometric points of $\X$. Let $\Phi_S \subseteq \PX$ be the normal
closure of the subgroupoid generated by images of $\ox$ for $x \in
S$. Then, the (isomorphism classes) of pointed covering spaces of
$\X$  that are FPR at all geometric points $y$ of $\Y$ lying above
points in $S$, are in one-to-one correspondence with open-closed
subgroups of $\pX$ that contain  $\Phi_S(x)$.
\end{prop}

This proposition will be used later in the proof of Theorem
\ref{T:finite}.
\section{Some Galois Theory of Gerbes}{\label{S:gerbes}}

 In this section we introduce a certain class of gerbes, called
 {\em monotonous gerbes} (Definition \ref{D:monotonous}), and
 study their Galois theory. The point is that the maps $\ox \:
 \phX \to \pX$ behave nicely for monotonous gerbes. We will see in
 the next section that every algebraic stack has a stratification
 by monotonous gerbes.

 We recall a fact from (\cite{SGA3}, Expos\'{e} VIB): Let $G \to
 X$ be a flat group space  of finite type. Let $G^o$ be the
 functor of {\it connected component of the identity}; that is,
 the functor that associates to any $X$-scheme $T$, the set of all
 $T$ points  $g \in G(T)$ of $G$ for which $g(t)$ falls within the
 connected component of identity of $G_t$ for every $t \in T$.
 This functor is a subsheaf of $G$. If this functor is
 representable, then $G^o$ is an open subgroup space of $G$.

 \begin{defn}{\label{D:monotonous}}
 We say that a finite type flat group space $ G \to X$ is
 {\em monotonous}, if the number of geometric connected components
 of the fiber of point $x \in X$ is a locally constant
 function on $X$, and if the functor $G^o$ defined above is representable
 (necessarily by an open subspace
 of $G$).  We say that an algebraic stack $\X$ is a {\sl monotonous gerbe},
 if it is connected and its stabilizer group is  monotonous.
 \end{defn}

It follows from (\cite{LM}, Corollaire 10.8), that a monotonous
gerbe is indeed an fppf gerbe, that it has a ``moduli space''
$\Xm$, and that the moduli map $\X \to \Xm$ is smooth and of
finite type. In fact, it is easy to see that this ``moduli space''
is a moduli space in the sense of Section\ref{S:moduli}, but we
will not need this fact here.

 If $G$ is a monotonous group space over $S$, then $G/G^o$ is a
 group space that is unramified and of constant degree over $S$.
 Therefore, it is  finite \'{e}tale over $S$. The following lemma
 is now immediate.

 \begin{lem}{\label{L:comp}}
 Let $X$ be a scheme, and let $G \to X$ be a  monotonous group space. Let $H$ be a subgroup
 that is both closed and open. Then, $H$ is monotonous.
 \end{lem}

 \begin{prop}{\label{P:cover}}
 Let $\X$ be a connected
 monotonous gerbe, and let $f \: \Y \to \X$ be a connected covering space. Then:
 \begin{itemize}
 \item[$\mathbf{i}$)] $\Y$ is a monotonous gerbe.
 \item[$\mathbf{ii}$)] If  $f$ is FPR (Definition \ref{D:FPR}) at some geometric point of $\Y$,
 then it is FPR everywhere.
 \end{itemize}
 \end{prop}

 \begin{proof}
 Let $p \: X \to \X$ be a chart for $\X$ and let $q  \: Y \to \Y$ be a connected component of
 the chart for $\Y$ obtained by pulling back
 $X$ via $f$. By Corollary \ref{P4.9}, the natural homomorphism $S_Y \to S_X \x_X Y$ is an open closed
 immersion. So it follows from Lemma \ref{L:comp}, $S_Y \to Y$ is monotonous. This proves ($\mathbf{i}$).
 Since $S_Y \to S_X \x_X Y$ is an open closed immersion of monotonous spaces, the necessary and
 sufficient condition for it to be an isomorphism is that the number of the geometric components of
 the fibers of $S_Y$ and $S_X \x_X Y$  be equal. But this can be checked at a single point.
 This proves ($\mathbf{ii}$).
 \end{proof}

 \begin{prop}{\label{P:normal}}
 Let $\X$ be a monotonous gerbe. Then, the image of $\OX \: \PhX \to \PX$ is a normal subgroupoid.
 In particular, the image of $\ox$ is a normal subgroup of $\pX$ for every
 geometric point $x$, and all these images are (non canonically) isomorphic.
 \end{prop}

 \begin{proof}
 Fix a point $x$ in $\X$. Let $I \subseteq \PX$ be the image of $\OX$, and let $N \subseteq \PX$ be its
 normal closure. Using Proposition \ref{P:b}, we can translate Proposition \ref{P:cover} as saying that,
 an open-closed subgroup of $\pX$ contains $N(x)$ if and only if it contains its subgroup $I(x)$.
 This implies that the closure of $N(x)$ and $I(x)$ are equal. But $I(x)$, being finite, is already closed.
 So $I(x)=N(x)$, which proves the claim.
 \end{proof}

\begin{cor}{\label{C:seq}}
Let $(\X,x)$ be a monotonous gerbe. Then we have the following
exact sequence:

$$\pi_1^{h}(\X,x) \to \pX \to \pXm \to 1.  $$

\end{cor}

\begin{proof}
This follows from Proposition \ref{P:normal} and  Theorem
\ref{G5}.
\end{proof}

\begin{rem}
The moduli map $\pi_{\X} \: \X \to \Xm$ of a monotonous gerbe
resembles a ``fibration with connected fibers'', and the above
exact sequence is the analogue of the ``fiber homotopy exact
sequence''.
\end{rem}

 The monotonous stacks for which the above sequence is short exact
 have quasi-finite stabilizers (Example \ref{P5.9}). If we
 restrict ourselves to reduced Deligne-Mumford stacks, then there
 is a simple description of the category of such stacks in
 (\cite{Noohi}, Section 6) which says that a stack with this
 property is uniquely determined by its moduli space, together
 with the above short exact sequence. In other words, to give a
 reduced monotonous Deligne-Mumford gerbe is the same as to give
 an algebraic space plus an extension of its fundamental group by
 a finite group. In particular, this gives a simple description of
 the category of zero-dimensional reduced Deligne-Mumford stacks
 (\cite{Noohi}, Section 6.1).

 \begin{cor}{\label{C:im}}
 Let $\X$ be an algebraic stack, and let $\X_0 \hookrightarrow \X$ be a locally closed immersion
 (or, simply, a monomorphism) of stacks. Assume $\X_0$ is monotonous.
 Let $x$ and $x'$ be geometric points of $\X$ that factor through $\X_0$.  Then, the image of
 $\ox \: \phX \to \pX$ is
 isomorphic (non canonically) to the image of $\omega_{x'} \: \pi_1^h(\X,x') \to \pi_1(\X,x')$.
 In particular, if $\ox$ is the zero map, then
 so is every other $\omega_{x'}$. Finally, if $\ox$ is injective, then so is $\omega_{x'}$.
 \end{cor}

\begin{proof}
Let $I_0 \subseteq \pi_1(\X_0,x)$ be the image of $\pi_1^h(\X_0,x)
\to \pi_1(\X_0,x)$, and $I \subseteq \pi_1(\X,x)$ be the image of
$\pi_1^h(\X,x) \to \pi_1(\X,x)$. Define $I'_0 \subseteq
\pi_1(\X_0,x')$ and $I' \subseteq \pi_1(\X,x')$ similarly. Note
that $\pi_1^h(\X_0,x)$ is canonically isomorphic to
$\pi_1^h(\X,x)$ (Proposition \ref{P:mono}). So $I$   is the image
of $I_0$   under the the natural map $\pi_1(\X_0,x) \to \pX$
(similarly for $I'$). Fix a ``path'' from to $x$ to $x'$ in
$\X_0$. We saw in the proof of the previous proposition, that this
``path'' induces an isomorphism $\pi_1(\X_0,x) \risom
\pi_1(\X_0,x')$, that maps $I_0$ isomorphically to $I'_0$. The
image of this ``path'' in $\X$, induces an isomorphism
$\pi_1(\X,x) \risom \pi_1(\X,x')$, that, because of the
commutativity of the following diagram, maps $I$ isomorphically
onto $I'$: $$\xymatrix{\pi_1(\X_0,x)   \ar[r] \ar[d]_{\sim} &
\pi_1(\X,x)  \ar[d]^{\sim} \\
            \pi_1(\X_0,x')  \ar[r]               &  \pi_1(\X,x')}$$
To prove the final statement, note that, if $\pi_1^h(\X,x)$ is
injective, then the stabilizer group of $x$ is zero dimensional
(Example \ref{P5.9}) and its set of geometric connected
components, which is a finite set, is isomorphic to $\phX$.
Therefore, since the stabilizer group $\mathcal{S}_{\X_0} \to
\X_0$  is flat, its relative dimension should be zero. In
particular, the stabilizer of $x'$ is also zero dimensional and
its set of geometric connected components is isomorphic to
$\pi_1^h(\X,x')$. Because $\mathcal{S}_{\X_0} \to \X_0$  is
monotonous, $\phX$ and $\pi_1^h(\X,x')$ have the same number of
elements. On the other hand, we just proved that their images
under, respectively $\ox$ and $\omega_{x'}$, have the same number
of elements. Therefore, if $\ox$ is injective, so will be
$\omega_{x'}$.
\end{proof}

 \section{Stratification by Gerbes and its Applications}{\label{S:strat}}

 The stratification theorem of this section is just a technical
modification of that of \cite{LM}, and so is its proof. We prove
that every algebraic stack has an stratification by monotonous
gerbes. We then use this stratification to prove a stronger
version of Theorem \ref{T:uniformization}.

 We will need the following lemma, due to Kai Behrend
 (\cite{Behrend}, Section 5), in the proof of our Stratification
 Theorem (Section \ref{S:strat}).

 \begin{lem}{\label{L:Behrend}}
  Let $G \to X$ be a  finite type group space, and assume $X$ is a Noetherian scheme.
  Then,  $X$ can be written as a disjoint union of  a finite
  family $(X_i)_{i \in I}$ of  reduced locally closed
 subschemes, such that for every $i \in I$, the restriction $f_i \: G \x_X {X_i} \to X_i$ is a group scheme
 for which the functor of {\em connected  components of the identity}
  (see the beginning of  Section \ref{S:gerbes}) is representable.
  \end{lem}

 \begin{lem}{\label{L:generic}}
 Let $f \: Y \to X$ be a morphism of finite presentation between Noetherian schemes. Then,
 $X$ can be written as a disjoint union of  a finite family $(X_i)_{i \in I}$ of  reduced locally closed
 subschemes such that for every $i \in I$, the restriction $f_i \: Y \x_X {X_i} \to X_i$ is  flat and has
 constant number of geometric fibers.
 \end{lem}

 \begin{proof}
 Lemma follows from the theorem on generic flatness (\cite{EGA4}, Th\'{e}or\`{e}me 6.9.1), and generic
 constancy of the number of geometric connected components (\cite{EGA4}, Proposition 9.7.8).
 \end{proof}

 \begin{thm}{\label{T:strat}}
 Let $\X$ be a Noetherian algebraic stack. Then, $\X$ can be written as a disjoint union of
 a finite family $(\X_i)_{i \in I}$ of locally closed subsets such that each $\X_i$ with its
 reduced structure is a monotonous gerbe.
 \end{thm}

 \begin{proof}
 In the proof given in (\cite{LM}, Th\'{e}or\`{e}me 11.5), simply substitute Lemma \ref{L:generic} and Lemma
 \ref{L:Behrend} above for the  `generic flatness theorem', and proof goes through verbatim.
 \end{proof}

Now we can prove the stronger version of Theorem
\ref{T:uniformization}.

\begin{thm}{\label{T:main}}
Let $\X$ be a Noetherian algebraic stack. Then, there exists a
covering space $\Y \to \X$ such that $\oy$ is the zero map for all
geometric points $y$ of $\Y$.
\end{thm}

\begin{proof}
We may assume $\X$ is connected. Let $(\X_i)_{i \in I}$ be the
stratification of $\X$ by monotonous gerbes (Theorem
\ref{T:strat}). Let $S$ be a finite subset of $|\X|$ that
intersects all $\X_i$. Fix a base point $x$ for $\X$. For any $x'
\in S$, choose a path connecting $x'$ to $x$, and use that to
transfer the image of $\omega_{x'}$ onto a  subgroup of $\pX$. The
union of all these subgroups is a finite subset of $\pX$ (see
Example \ref{P5.9}), so we could choose an open normal subgroup of
finite index $U \subseteq \pX$ that does not intersect any of
these subgroups. The corresponding covering space will have the
desired property (Lemma \ref{P6} and Corollary \ref{C:im}).
\end{proof}

 Let $S$ be a subset of $|\X|$ that intersects all the $\X_i$.
 Assume $\ox$ is the zero map for all $x$ in $S$. As in the proof
 of the above theorem, Lemma \ref{P6} together with Corollary
 \ref{C:im} imply that $\ox$ is the zero map for all $x$ in $\X$.
 In particular, let $S \subset |\X|$ be a set of points in $\X$
 that has the following property: For every locally closed subset
 $T$ of $|\X|$ , the intersection $T \cap S$ is nonempty. Then,
 once the injectivity of $\ox$ is established for all $x \in S$,
 then it follows that $\ox$ is injective for all $x$. For
 instance, the set of closed point of an algebraic stack of finite
 type over a field has this property.

Let $\Y \to \X$ be a covering space of $\X$. Let $y$ be a point in
$\Y$ and $x=f(y)$ its image in $\X$. Theorem \ref{T:main} is most
useful when combined with the following facts (that are always
true):

\begin{itemize}

\item  The kernel of $\oy$ is the same as the kernel of $\ox$.

\item The stabilizer of $\Y$ is an open closed subgroup of the
pull back of the stabilizer of $\X$.

\end{itemize}

 For example, if in Theorem \ref{T:main} we assume in addition
 that $\ox$ is injective for every geometric point $x$ of $\X$ ,
 then it follows that the hidden fundamental groups of $\Y$ are
 all trivial. In other words, the stabilizer group of $\Y$ has
 zero-dimensional geometrically connected fibers. If we further
 assume that $\X$ has unramified stabilizer (i.e., $\X$ is
 Deligne-Mumford), then so does $\Y$. Therefore, the stabilizer of
 $\Y$ will be trivial, which means that $\Y$ is an algebraic space
 (compare Theorem \ref{T:uniformization}).

\section{A Finiteness Theorem}

All the algebraic stack in this section are assumed to be
Noetherian.  We prove the following

\begin{thm}{\label{T:finite}}
Let $(\X,x)$ be a connected algebraic stack that has a moduli
space (in the sense of Section \ref{S:moduli}). Let $N$ be the
kernel of the $\pX \to \pXm$.  Then, there exists a finite set $T$
of torsion elements in $N$ such that $N$ is the smallest closed
normal subgroup containing $T$.
\end{thm}

 The main ingredient of the proof is the following result, which
 can be thought of as a generalization of Proposition
 \ref{P:cover}

\begin{prop}{\label{P:cover2}}
Let $\X$ be and algebraic stack Let $f \: \Y \to \X$  be a Galois
covering space. Let $(\X_i)_{i \in I}$ be a stratification of $\X$
by monotonous substacks (Theorem \ref{T:strat}). Let $S$ be a
subset of  $\X$ that intersects all the $\X_i$. Assume the maps
$\phY \to \phX$ are isomorphisms for all geometric points $y$ of
$\Y$ lying above points in $S$. Then,  $f$ is fixed point
reflecting.
\end{prop}

\begin{proof}
Let $\Y_i =\X_i \x_{\X} \Y$. Clearly $\Y$ is a disjoint union of
$\Y_i$. Since the hidden fundamental group of a point only depends
on its residue gerbe, it is enough to prove that  each $\Y_i \to
\X_i$ is fixed point reflecting. Each $(\Y_i)$ is a disjoint union
of connected algebraic stacks that are all isomorphic to each
other (because the action of the Galois group permutes them around
transitively). So it is enough to prove the proposition for one of
these connected components. Such a connected component is a
monotonous gerbe by Proposition \ref{P:cover} ($\mathbf{i}$), so
the result follows from Proposition \ref{P:cover} ($\mathbf{ii}$)
\end{proof}

\begin{proof}[Proof of Theorem \ref{T:finite}]
Fix a finite set $S$ as in Proposition \ref{P:cover2}. In this
proof, we will switch our notation for elements of $S$ from $x$ to
$s$. Fix a geometric point $x$ for $\X$ all through the proof. Let
$\Phi \subseteq \PX$ be the closed normal subgroupoid of $\PX$
generated by the images of $\ox$ for various $s \in S$, and let
$N=\Phi(x)$. The group $N$ has the desired finiteness property of
the theorem. We claim that $\pX/{N}$ is naturally isomorphic to
$\pXm$. The set of open-closed normal subgroups of $\pXm$ is in
one-to-one correspondence with the set of Galois covering spaces
of $\Xm$. By Proposition \ref{G9}, the set of Galois covering
spaces of $\Xm$ is in one-to-one correspondence with the set of
fixed point reflecting Galois covering spaces of $\X$. By
Proposition \ref{P:cover2} and  Proposition \ref{P:b}, this set is
in one-to-one correspondence with the set of open-closed normal
subgroupoids of $\PX$ that contain the image of $\ox$ for all $s
\in S$ (equivalently, contain $\Phi$). This set is in one-to-one
correspondence with open-closed normal subgroups of $\pX$
containing $N$. Claim now follows from Lemma \ref{G10.5}.
\end{proof}

 \pagebreak


\setcounter{section}{1} \setcounter{equation}{1}
\renewcommand{\thesection}{\Alph{section}}

\section*{Appendix: Fibrations}


In this section we introduce a notion of a {\em fibration} between
algebraic spaces and prove that being a fibration is local in the
fppf topology. The results of this appendix are used in section
\ref{S:application}.  This appendix is divided into four
subsections. In the first part, we define fibrations in terms of
the fiber homotopy exact sequence of fundamental groups. In the
second part, we discuss the geometric meaning of being a
fibration. In the third  part, we prove that fibrations satisfy
fppf descent. In the last part we discuss
 fibrations  of algebraic stacks.

In this appendix all algebraic spaces and morphisms between them
are assumed to be quasi-compact.

\subsection{Fibrations: the definition}

 Let $f \: Y \to X$ be a morphism of algebraic spaces.
We would like to study those $f$ for which  the following axiom
holds:

 \begin{itemize}
\item[$\mathbf{A}$)]   For any choice of base point $y$ and $x=f(y)$, the sequence
\begin{equation}{\label{3}}
\pi_1(Y_x,y) \to \ppYY \to \ppXX \to \pi_0(Y_x,y) \to \pi_0(Y,y)
\to \pi_0(X,x) \to \{*\}
\end{equation}
 is exact.
 \end{itemize}

 When $X$ and $Y$ are connected, the above exact sequence takes the following form:
 \begin{equation}{\label{2}}
 \pi_1(Y_x,y) \to \ppYY \to \ppXX \to \pi_0(Y_x,y) \to \{*\}
 \end{equation}

The axiom $\mathbf{A}$ means that this shorter sequence is exact
for the restriction of $f$ to any connected component of $X$ (and
the corresponding component of $Y$).

\begin{defn}{\label{D:fibration}}
 A surjective morphism $f \: Y \to X$  is called a {\em quasi-1-fibration},
 or simply a {\sl fibration}, if every base extension of $f$ satisfies $\mathbf{A}$.
\end{defn}

A well-known class of  fibrations are covering (i.e., finite
\'{e}tale) morphisms. More generally, it is proven in
(\cite{SGA1}, Expos\'{e} X, 1.6) that any proper flat morphism $f
\: Y \to X$
 that has geometrically reduced fibers  is a fibration ($X$ locally Noetherian).
The following proposition is a rather tedious exercise in diagram
chasing:

\begin{prop}{\label{P:comp}}
Let $g \: Z \to Y$ be a fibration, and let $f \: Y \to X$ be
morphisms that satisfies $\mathbf{A}$.  Then $f \circ g$ satisfies
$\mathbf{A}$. If $f$ and $g$  are both fibrations, then so is $f
\circ g$.

\end{prop}

\subsection{Geometric interpretation}

\begin{defn}
By a {\em geometrically connected} morphism of algebraic spaces we
mean the one whose geometric fibers are connected.
\end{defn}

In this section we analyze the geometric meaning of $\mathbf{A}$
for geometrically connected morphisms $f \: Y \to X$. In this
case, the exact sequence (\ref{3}) takes the following form
 \begin{equation}{\label{1}}
\pi_1(Y_x,y) \to \ppYY \to \ppXX \to 1.
\end{equation}

Let us   translate this exactness into geometry.  Assume $f \:
(Y,y) \to (X,x)$ is a pointed morphism (not necessarily
geometrically connected) of connected algebraic spaces.  Then we
have the following

\begin{prop}{\label{P:surj}} Let $X'$ be  covering space of $X$,
and let $U \subset \ppXX$ be the corresponding subgroup (for some
choice of base point). Let $Y'$ be the pull-back of $X'$ via $f$.
Then, there is a natural bijection between the set of connected
components of $Y'$ and the set of orbits of the (right) action
(along $f$) of $\ppYY$ on the set $U\backslash\ppXX$ of right
cosets of $U$. In particualr, $\ppYY \to \ppXX$ is surjective if
and only if the pull back (via $f$) of  every (connected)
 covering space  $X'$ of $X$ is connected.
\end{prop}

\begin{prop}{\label{P:first}}
The sequence (\ref{1}) being exact is equivalent to the following:
\begin{itemize}
\item A covering space $P: Y' \to Y$ descends to a  covering space
of $X$ if and only if the restriction of $Y'$ over the  geometric
fiber $Y_x$ has a section. Furthermore, the descended covering
space is unique (up to isomorphism).
\end{itemize}
\end{prop}

 It is easy to see that the uniqueness
of descent in ($\mathbf{ii}$) is valid for non connected finite
\'{e}tale covers as well.

Let us recall a well-known

\begin{lem}{\label{L:a}}
Let $f \: (Y,y) \to (X,x)$ be a pointed morphism of connected
algebraic spaces. Let $(X',x') \to (X,x)$ be a pointed covering
space. Then, $f$ lifts (necessarily uniquely) to $(X',x')$ if and
only if the image of $\ppYY$ in $\ppXX$ under $f$ is contained in
the image of $\pXX$.
\end{lem}

\begin{proof}
One implication is trivial.  To prove the non trivial one, assume
the image of $\ppYY$ in $\ppXX$ under $f$ is contained in the
image of $\pXX$. Consider the pull back of $(X',x')$ along $f$.
This pull back is a disjoint union of covering spaces of $(Y,y)$,
and it is naturally pointed. Let $(Y',y')$ be the connected
component  that contains the base point. It corresponds to the
open subgroup of $\ppYY$ which, by hypothesis, is equal to the
entire $\ppYY$. So $(Y',y') \to (Y,y)$ is trivial. The inverse map
$(Y,y) \to (Y',y')$ composed with $(Y',y') \to (X',x')$ is the
required lift.
\end{proof}

\begin{prop}{\label{P:decomp}}
Let  $f \: Y \to X$ be a morphism of connected algebraic spaces
that satisfies $\mathbf{A}$. Then, there is a factorization $f =p
\circ g$ such that $p$ is finite \'{e}tale and $g$ is
geometrically connected. Furthermore, $g$ satisfies $\mathbf{A}$.
The decomposition $f =p \circ g$ is unique up to a unique
isomorphism and commutes with base change.
\end{prop}

\begin{proof}
For any choice of a base point $y$ for $Y$ we have the following
exact sequence: $$\pi_1(Y_x,y) \to \ppYY \to \ppXX \to
\pi_0(Y_x,y) \to \{*\},$$ where $x=f(y)$. Since the map $\ppYY \to
\ppXX$ is independent of the choice of the base-point $y$ (up to
isomorphism), the cardinality of $\pi_0(Y_x,y)$ will also be
independent of $y$. So if we move $y$ around, we see that, as long
as $x$ is in the image of $f$, the number of geometric connected
components of $Y_x$ remains constant (equal to
$\#|\pi_0(Y_x,y)|$). Now fix the base-point $y$. Let $U$ be the
image of $\ppYY$ in $\ppXX$, and let $p \: (X',x') \to (X,x)$ be
the covering space associated to it. This covering has degree
$\#|\pi_0(Y_x,y)|$, and, by Lemma \ref{L:a}, there is a unique
factorization $f =p \circ g$. Clearly, $g$ is geometrically
connected. It is also easy to check that, any other factorization
$f =p' \circ g'$ in which $p'$ is finite \'{e}tale and $g$
geometrically connected is uniquely isomorphic to the
factorization $f =p \circ g$ (to see this, consider the open
subgroup of $\ppXX$ corresponding to $p'$. It should be contained
in $U$, and should have the same index; so it must be equal to
$U$). We have proved so far the existence and uniqueness of the
factorization. We now prove that $g$ satisfies $\mathbf{A}$. By
the very construction of $g$, the axiom $\mathbf{A}$ is satisfied
at the specific base-point $y$. But, by uniqueness, changing $y$
does not change the factorization (up to isomorphism); so the
exactness of the above sequence is valid for any other choice of
the base point as well. \\ That the factorization commutes with
base change is trivial.
\end{proof}

\subsection{Fibrations and Descent}

Fibrations behave well with respect to fppf descent:

\begin{thm}{\label{T:descent}}
Let $f \: Y \to X$ be a morphism of algebraic spaces. Assume $a \:
V \to X$ is a (surjective) morphism that is universally effective
descent morphism for finite \'{e}tale morphisms. Further, assume
$a$ is either universally open or universally closed. If the base
extension of $f$ along $a$  is a fibration, then so is $f$.
(Recall that  `universally'  means that every base extension of
$a$ is also an effective descent morphism.)
\end{thm}

\begin{proof}
We may assume $X$ and $Y$ are connected. Let $W=Y \x_X V$, and let
$r \: W \to V$ be the projection map. $$\xymatrix{ W \ar[r]^{b}
\ar[d]_{g} & Y \ar[d]^{f} \\
             V \ar[r]_{a}        & X }$$

By Proposition \ref{P:decomp}, there is a unique (up to a unique
isomorphism) factorization $g = q \circ t$ such that $q$ is finite
\'{e}tale, and $t$ has geometrically connected fibers. The
uniqueness of this factorization, plus the fact that it is
invariant under base change, imply that it satisfies descent
condition relative to the map $a$ (namely, by repeating the same
argument with $V \x_X V$ and $V \x_X V \x_X V$ instead of $V$).
Therefore, since $a$ is a universally effective descent morphism,
the factorization descends to a factorization for $f$. Hence, we
are reduced to the case where $f$ has geometrically connected
fibers. In this situation, the exact sequence (\ref{2}) implies
that $g \: W \to V$ is a disjoint union of geometrically connected
morphisms that satisfy $\mathbf{A}$.

To prove the exactness of (\ref{1}) for $f$, we can now use the
geometric translation of the exactness as in Proposition
\ref{P:first}.
 Let $Y' \to Y$ be a (connected) covering space of $Y$. Pick a
 geometric point $x$ for $X$, and assume $Y' \to Y$
 has a section over $Y_x$. We want to show that $Y'$ descends
 uniquely to a covering space $X'$ of $X$.
Pick a connected component $V_0$ of $V$ whose image in $X$
contains $x$, and let $W_0$ be the (unique) connected component of
$W$ above it. Pick a geometric point $v$ in $V_0$ lying above $x$.
Let $W'=W \x_Y Y'$  be the pull back of $Y'$  to $W$, and let
$W'_0$ be its restriction to $W_0$. Then, the restrictions of
$W'_0$ to $W_{0,v}$ has a section. Hence, $W'_0$   descends
(uniquely) to a covering space  $V'_0$
  of $V_0$. This, in particular, implies that, if we replace $x$
  by any other point in the image of $V_0$, then $Y' \to Y$ has a
  section over $Y_x$. Therefore, if $V_1$ is another connected
  component of $V$ whose image in $X$ intersects that of $V_0$,
  then $W'_1$ will descend (uniquely) to a covering space  $V'_1$
  of $V_1$. We can continue this argument inductively and, using
  the fact that $X$ is connected and the map $V \to X$ is open (or
  closed), we see that any point in $X$ can be eventually reached
  after finitely many steps. Therefore, we have shown that $W'$
  descends uniquely to a covering space  $V'$ of $V$.

  Repeating the same argument  with $V \x_X V$ and $V \x_X V \x_X V$
  instead of $V$, and using the uniqueness
of descent, we see that $V'$   satisfy the descent condition
relative to $a$, so it descends to to a covering space $X'$ of
$X$. We claim that $X'$ does the job. Let $Y''$ be the pull back
of $X'$ to $Y$. We claim that $Y''$ is isomorphic to $Y'$. To do
so, it is enough to show that their pull backs to $W$ are
isomorphic (because, $b \: W \to Y$ is an effective descent
morphism). But this is true by the construction of $X'$. It only
remains to prove that descent is unique; namely, that if $X'$ and
$X''$ are two covers of $X$ whose pull backs $Y'$ and $Y''$ to $Y$
are isomorphic, then they are isomorphic. It is enough to show
that pull backs $V'$ and $V''$ of
 $X'$ and $X''$ to $V$ are isomorphic. But, since $g \: W \to V$ is
 a disjoint union of geometrically
 connected morphisms that satisfy $\mathbf{A}$, it is enough to prove that  the pull backs
 $W'$ and $W''$ of $V'$ and $V''$ to $W$ are isomorphic (Proposition \ref{P:first}).
 But, these two are isomorphic to pull backs of $Y'$ and $Y''$, which are already
 isomorphic over $Y$, to $W$.
 The proof is now complete.
\end{proof}

 The following classes of morphisms satisfy the requirements of the Theorem
 (\cite{SGA1}, Expos\'{e} IX):

\begin{itemize}
\item fppf.

\item Proper, surjective and of finite presentation.

\item Surjective, universally open and of finite type (target locally Noetherian).
\end{itemize}

\subsection{Fibrations for algebraic stacks}{\label{S:3}}

The definition of a fibration and the arguments of the previous
sections are quite formal and can be generalized algebraic stacks.
\begin{defn}{\label{D:stfib}}
We say that a  morphism $f \:\Y \to \X$ of algebraic stacks is a
{\em fibration} if for any choice of base points, the sequence
(\ref{3}) is exact.
\end{defn}

In this definition choice of base points $y$ and $x$ should be
made so that $f(y)$ is 2-isomorphic to $x$. A clarification here
is in order: A pointed morphism between pointed stacks $(\Y,y)$
and $(\X,x)$ consists of a morphism $f \:\Y \to \X$ together with
a transformation $\phi \: x \rsa f(y)$. We do need to have the
extra data $\phi$ in order to be able to form the sequence
(\ref{3}). However, the exactness of this sequence is independent
of the choice of $\phi$.

When $f$ is a representable morphism, then $f$ being a fibration
is equivalent to its base extension via every $X \to \X$ being a
fibration, where $X$ is a scheme (or an algebraic space). In fact,
by Theorem \ref{T:descent}, to check that $f$ is a fibration we
can pick a flat chart $X \to \X$ and check whether the base
extension of $f$ to $X$ is a fibration. This implies the following

\begin{thm}{\label{T:last}}
Let $R \st{s,t} X$ be a flat and locally of finite presentation
 groupoid so that $R \to X \x X$ is separated and quasi compact
 (these conditions ensure that the
 quotient stack is algebraic).
 Let $\X=[X/R]$ be the corresponding quotient stack. Then, for $X \to \X$ to be a fibration it
 is necessary and sufficient that $s$ (and $t$) be fibrations. In particular, let $G \to S$ be a  group space
 that is flat, separated, locally of finite presentation and a fibration over $S$.
 Assume $G$ acts on an algebraic space $X$ over $S$. Then, the quotient map $X \to [X/G]$ is a fibration.
 \end{thm}

\section*{Acknowledgments}
Part one of this paper is a modification of some part of my MIT
thesis. I would like to express my gratitude to my thesis advisors
S. Kleiman and F. Diamond. I learned algebraic stacks in a truly
fantastic course taught by J. de Jong at MIT in spring 1999. I am
forever indebted to him. The second part of the this paper was
added during my visit to UBC. I would like to thank Kai Behrend
for the invitation and for sharing many of his ideas.

\bibliographystyle{amsplain}
\bibliography{FundGp}
\end{document}